\documentclass[12pt,reqno]{amsart}
\usepackage[english]{babel}
\usepackage[T1]{fontenc}
\usepackage{amsthm,amsfonts,amssymb,amscd,color}
\usepackage[all]{xy}

\newtheorem{theorem}[subsection]{Theorem}

\newtheorem{proposition}[subsection]{Proposition}

\newtheorem{conjecture}[subsection]{Conjecture}
\newtheorem{lemma}[subsection]{Lemma}
\newtheorem{corollary}[subsection]{Corollary}

\theoremstyle{definition}
\newtheorem{definition}[subsection]{Definition}
\newtheorem{example}[subsection]{Example}
\newtheorem{proposition-definition}[subsection]{Proposition-Definition}

\theoremstyle{remark}
\newtheorem{remark}[subsection]{Remark}

\newcommand{\fatdot}{{\raisebox{1.pt}{$\scriptscriptstyle \bullet$}}}

\newcommand{\dual}{{{\scriptscriptstyle \vee}}}

\newcommand{\At}{\operatorname{At}\nolimits}
\newcommand{\ch}{\operatorname{ch}}

 \newcommand{\Coh}{\operatorname{Coh}\nolimits}

\newcommand{\EXT}{{\:\mathcal Ext\:}}
\newcommand{\Ext}{\operatorname{Ext}\nolimits}
\newcommand{\GL}{{GL}_5(\CC)}
\newcommand{\gr}{\operatorname{gr}}
\newcommand{\Gr}{\operatorname{Gr}}
\newcommand\Hilb{{\operatorname{Hilb}\nolimits}}

\newcommand{\HOM}{{{\mathcal H}om\:}}
\newcommand{\Hom}{\operatorname{Hom}\nolimits}

\newcommand{\im}{\operatorname{im}\nolimits}

\newcommand{\Pf}{\operatorname{Pf}\nolimits}

\newcommand{\pr}{\operatorname{pr}\nolimits}

\newcommand{\rk}{\operatorname{rk}\nolimits}
\newcommand{\Sing}{\operatorname{Sing}}

\newcommand{\TOR}{{\:\mathcal T\!or}}
\newcommand{\Tr}{\operatorname{Tr}}

\newcommand{\CC}{{\mathbb C}}
\newcommand{\HH}{{\mathbb H}}

\newcommand{\ZZ}{{\mathbb Z}}
\newcommand{\QQ}{{\mathbb Q}}
\newcommand{\PP}{{\mathbb P}}

    \newcommand{\LL}{{\mathbb L}}

\newcommand{\OOO}{{\mathcal O}}
    \newcommand{\CO}{{\mathcal O}}

    \newcommand{\CI}{{\mathcal I}}
\newcommand{\GGG}{{\mathcal G}}
    \newcommand{\CG}{{\mathcal G}}
\newcommand{\EEE}{{\mathcal E}}
    \newcommand{\CE}{{\mathcal E}}
    \newcommand{\D}{{\mathcal D}}

    \newcommand{\CH}{{\mathcal H}}

\newcommand{\FFF}{{\mathcal F}}

    \newcommand{\CF}{{\mathcal F}}
\newcommand{\CCC}{{\mathcal C}}

\newcommand{\NNN}{{\mathcal N}}
    \newcommand{\CN}{{\mathcal N}}
\newcommand{\MMM}{{\mathcal M}}

\newcommand{\RRR}{{\mathcal R}}

\newcommand{\TTT}{{\mathcal T}}
    \newcommand{\CT}{{\mathcal T}}

\newcommand{\FFFF}{{\boldsymbol{\mathcal F}}}

\newcommand\alp{\alpha}

\newcommand\eps{\epsilon}

\newcommand\IFF{{\ \Longleftrightarrow\ }}
\newcommand{\into}{\hookrightarrow}

\newlength{\rrrr}
\newcommand{\isom}[1]{{\settowidth{\rrrr}{$\scriptstyle{x#1x}$}
\xrightarrow{\makebox[\rrrr]{$\scriptstyle{#1}$}}
\hspace{-0.5\rrrr }\hspace{-1.1 em}
\raisebox{- 0.5 ex}{$\sim$}\hspace{0.7\rrrr }
}}

\newcommand\lra{{\longrightarrow}}

\newcommand\rar{\rightarrow}
\newcommand\qis{{\:\raisebox{-1.2 ex}
{$\stackrel{\textstyle\simeq}{\scriptstyle\mathrm{qis}}$}\:}}

\newcommand{\simvert}{\mbox{\raisebox{0.8 ex}
{{\scriptsize\rm )}\hspace{-1.ex}\raisebox{.7em}{\scriptsize\rm (}}\rule{0pt}{1pt}}}

\newcommand{\tr}{\mathop{{\mathsf{tr}}}}

\renewcommand\emptyset{\varnothing}

\renewcommand{\bar}[1]{\overline{#1}}

    \newcommand{\AJ}{{{\mathrm{AJ}}}}
    
    \newcommand{\td}{{\mathop{\mathrm{td}}}}
    \newcommand{\TY}{{\tilde{Y}}}

\newcommand{\CS}{{\mathcal{S}}}
\newcommand{\CA}{{\mathcal{A}}}
\newcommand{\CB}{{\mathcal{B}}}

\newcommand{\CM}{{\mathcal{M}}}
\newcommand{\TF}{{\tilde{F}}}
\newcommand{\TTH}{{\tilde{H}}}

\newcommand{\Ker}{\operatorname{Ker}\nolimits}
\newcommand{\rank}{\operatorname{rank}\nolimits}
\renewcommand{\GL}{\operatorname{GL}\nolimits}
\newcommand{\CY}{{\mathcal{Y}}}
\newcommand{\CV}{{\mathcal{V}}}
\newcommand{\HOH}{{\mathsf{HH}}}
\newcommand{\HI}{{\widehat{I}}}
\newcommand{\MI}{{\mathsf{MI}}}

\author{A. Kuznetsov}

\address{A. K.: Algebra Section, Steklov Mathematical Institute, 8 Gubkin str., Moscow 119991 Russia}
\email{akuznet@mi.ras.ru}

\thanks{A.K. was partially supported by
RFFI grants 05-01-01034, 07-01-00051, 07-01-92211, and NSh-9969.2006.1,
INTAS 05-1000008-8118,
the Russian Science Support Foundation,
and gratefully acknowledge of the support of the Pierre Deligne fund
based on his 2004 Balzan prize in mathematics.}

\author{L. Manivel}

\address{L. M.: Institut Fourier, UMR 5582 UJF-CNRS, Universit\'e Joseph
Fourier, F-38402 Saint Martin d'H\`eres cedex, France}
\email{Laurent.Manivel@ujf-grenoble.fr}

\author{D. Markushevich}

\address{D. M.: Math\'ematiques - b\^{a}t. M2, Universit\'e Lille 1,
F-59655 Villeneuve d'Ascq Cedex, France}
\email{markushe@math.univ-lille1.fr}


\title{Abel-Jacobi maps for hypersurfaces \\ and non commutative Calabi-Yau's}

\begin{document}
\begin{abstract}
It is well known that the Fano scheme of lines on a cubic 4-fold is a symplectic variety.
We generalize this fact by constructing a closed $(2n-4)$-form on the Fano scheme of lines
on a $(2n-2)$-dimensional hypersurface $Y_n$ of degree $n$. We provide several definitions
of this form --- via the Abel--Jacobi map, via Hochschild homology, and via the linkage
class --- and compute it explicitly for $n = 4$.
In the special case of a Pfaffian hypersurface $Y_n$ we show that the Fano scheme is birational
to a certain moduli space of sheaves of a $(2n-4)$-dimensional Calabi--Yau variety $X$
arising naturally in the context of homological projective duality, and that the constructed form is induced by the holomorphic volume form on $X$. This remains true for a general non Pfaffian hypersurface but the dual Calabi-Yau becomes non commutative.
\end{abstract}
\maketitle

\section*{Introduction}

Let $Y$ be a projective or compact K\"ahler manifold
and $p:Z\rar B$ a family of $k$-cycles on $Y$, parameterized by a smooth base
$B$. The Abel--Jacobi map (or the cylinder map) of the family $p:Z\rar B$ is the homomorphism
$$
\AJ : H^\fatdot(Y,\ZZ)\lra H^{\fatdot-2k}(B,\ZZ)\  ,\ \ \ \AJ(c)=p_*q^*(c),
$$
where $q:Z\rar Y$ is the natural projection.

If $Y$ is a nonsingular hypersurface in $\PP^{n}$, then the only
interesting piece of the cohomology of $Y$ is the primitive part
of $H^{n-1}(Y)$.  Clemens and Griffiths \cite{CG} studied the Abel--Jacobi map
when $Y$ is a cubic threefold in $\PP^4$ and $B=F(Y)$ is the Fano surface of
$Y$, that is the base of the family of lines on it. They showed that $\AJ$ is an isomorphism
between the Hodge structures on $H^3(Y)$ and $H^1(F(Y))$, and deduced
the nonrationality of $Y$. Beauville and Donagi \cite{BD} considered
the case of a smooth cubic fourfold in $\PP^5$. They proved that $\AJ$
provides an isomorphism of polarized Hodge structures between the primitive
cohomologies $H^4(Y)_{prim}$ and $H^2(F(Y))_{prim}$. In particular,
$H^{3,1}(Y)\simeq H^{2,0}(F(Y))$ is 1-dimensional, that is $F(Y)$ carries a holomorphic
2-form $\alp\in H^{2,0}(F(Y))$, unique up to proportionality. Looking at a special cubic $Y$,
whose equation is the Pfaffian of a 6-by-6 matrix of linear forms, they identified
$F(Y)$ with the Hilbert square $X^{[2]}$ of the ``orthogonal'' K3 surface $X=Y^\perp$,
which is a transversal linear section of $G(2,6)$.
They deduced from this that $\alp$ is nondegenerate, hence symplectic, and moreover,
that for any smooth cubic $Y$,
the Fano scheme $F(Y)$ is an irreducible symplectic fourfold, obtained by
deformation of the complex structure on $X^{[2]}$.

In \cite{KM}, we take the point of view
that though the ``orthogonal'' K3 surface $X$ is undefined for a general smooth
cubic $Y$, its derived category still makes sense. Following Bondal, we refer
to this ``derived category'' as a noncommutative (or categorical) K3 surface.
It was first identified in the paper \cite{Ku1} as the orthogonal complement $\CCC(Y)$
to the natural exceptional collection $\OOO_Y$, $\OOO_Y(1)$, $\OOO_Y(2)$ in the derived category $\D^b(Y)$. For brevity, we call the relationship between $Y$ and $X$ (or $\CCC(Y)$) described in \cite{BD,KM} by Beauville-Donagi correspondence.
The Fano scheme $F(Y)$ parametrizes, at the same time: lines in $Y$,
length-2 zero-dimensional subschemes of $X$, and certain objects in $\CCC(Y)$.
Thus it represents an important ingredient of the Beauville-Donagi correspondence.

The main result of the paper is a generalization of
the Beauville-Donagi correspondence to higher dimensions, where
the cubic is replaced by a hypersurface $Y_n$ of degree $n$ in $\PP^{2n-1}$,
either Pfaffian or general, and the K3 surface becomes a (noncommutative) Calabi-Yau.
The central role belongs to the Fano scheme $F(Y_n)$ parametrizing lines in $Y_n$:
it is the main object of applications and the main testing ground for our
techniques.
Several variants of Abel-Jacobi maps appear as
cohomological descents of Fourier-Mukai transforms. As soon as we want
to work with ideal ``varieties'' which have no points but have
derived categories, we use the Hochschild homology, well
defined both for varieties and for a certain class of triangulated categories.
Thus a part of the paper describes categorical tools, such
as universal Atiyah class, Hochschild homology and its functorial properties,
Hochschild--Kostant--Rosenberg (HKR) isomorphisms, Calabi-Yau categories
and Homological Projective Duality.
Another substantial part is a background material
about the Pfaffian varieties, their resolutions of singularities and
linear subspaces on them, the facts which we need and which
are either new, or are not easily extracted from the existing literature.

We will now enumerate the higher-dimensional analogs of the respective features
of the Beauville-Donagi correspondence which constitute our generalization.

{\em We are first assuming $Y_n$ general.}

1. The topological cylinder map $AJ$ of the
family of lines induces a Hodge isometry
$H^{2n-2}(Y_n,\ZZ)_{prim}\simeq H^{2n-4}(F(Y_n),\ZZ)_{van}/(\mathrm{torsion})$,
where the superscript $van$ denotes the vanishing cohomology, defined as the
kernel of the Gysin map associated to the embedding $F(Y_n)\hookrightarrow G(2,2n)$
(Proposition \ref{Hodge_isometry}).
Contrary to the Beauville-Donagi case ($n=3$), we do not know whether
the cohomology of $F(Y_n)$ is torsion-free, and the vanishing cohomology
may be strictly smaller than the primitive one: we prove that for $n=4$ it has
corank 1 (Proposition~\ref{Hodge_num}).

2. There is a noncommutative $(2n-4)$-Calabi-Yau $\CCC_n=\CCC(Y_n)$
naturally associated to $Y_n$ (Theorem \ref{cccn}), and the Fano
scheme $F(Y_n)$ is identified with a ``fine moduli space'' of
objects in $\CCC_n$ (since a general definition of a moduli space
in a triangulated category is still missing, we explain in Remark~\ref{moduli_complexes}
which features of $F(Y_n)$ enable us to consider it in this sense,
see also Proposition \ref{cgcc} where the objects of $\CCC_n$ parameterized
by points of $F(Y_n)$ are described). In the Beauville-Donagi case, $\CCC_3$ is a
deformation of the derived category of a K3 surface in the sense
that it becomes $\D^b(K3)$ for special cubics $Y$.

3. The higher-dimensional analog of Beauville-Donagi's symplectic form
is a closed exterior $p$-form $\alp_p$ with $p=2n-4$. It is obtained in several manners:
(a) as a generator of the \mbox{1-dimensional} vector space
$H^{2n-4,0}(F(Y_n))=\AJ(H^{2n-3,1}(Y_n))$ (Corollary \ref{2n-4});
(b) via the map on the Hochschild homology induced by the universal line
$Z\subset Y_n\times F(Y_n)$, followed by the HKR isomorphism (Corollary
\ref{hoch_atiyah});
(c) as the $p$-form induced on $F(Y_n)$ by the holomorphic
volume element of the noncommutative Calabi-Yau $\CCC_n$ (Proposition \ref{induced_cy});
(d) via the composition of Yoneda coupling with the Atiyah class or with the
linkage class (Theorem~\ref{epsat}, Corollary \ref{closed}).
Note that the approaches in (b) and (d) apply in a more general situation
and provide a $p$-form on the Hilbert scheme of $Y_n$ and on moduli spaces
of sheaves on $Y_n$.

4. The nondegeneracy of Beauville-Donagi's symplectic form has no
obvious analog in higher dimension. We found two natural ways
to measure the nondegeneracy of a $p$-form $\alp$ on a \mbox{$N$-dimensional} vector space $V$ ($N\geq p$) with even $p=2k$: (a) the 2-rank $r^{(2)}$ of $\alp$, that is the rank of the induced bilinear form on $\wedge^k V$;
(b) the dimension $d$ of the orbit of $\alp$ in $\wedge^pV$ under the action
of $GL(V)$. For particular values of $N$, $k$, there are alternative notions, for example:
(c) if $N=7$ and $k=4$,
there is a natural way to make a symmetric bilinear form~$q_\alp$ on~$V$
out of $\alp$ (see Section \ref{4-forms}),
and the rank $r$ of $\alp$ is defined to be the rank of~$q_\alp$.
One says that
$\alp$ is nondegenerate if $r^{(2)}$ (resp. $d$, or
$r$) is maximal for given values of $N, k$.
We explicitly compute $\alp_4$ on the 7-fold $F(Y_4)$ for a 6-dimensional quartic
$Y_4$ (Theorem \ref{determinants},
Corollary \ref{4-form}) and determine the three invariants for it:
$(r^{(2)},d,r)=(18,34,4)$, whilst the maximal values are $(21,35,7)$, corresponding to
the 4-forms which fill an open orbit and whose stabilizer is the exceptional linear groupe $\mathbf G_2$ (Section \ref{4-forms}). Our result can be interpreted by saying that $\alp_4$ is minimally degenerate, for it belongs to a codimension-1 orbit of $GL(7)$.

{\em Let now $Y_n$ be a generic Pfaffian hypersurface in} $\PP^{2n-1}$, that is the zero locus of the Pfaffian of a generic $(2n)\times(2n)$ skew-symmetric matrix of linear
forms in $2n$ variables.

5. We associate to $Y_n$ an ``orthogonal'' Calabi--Yau manifold $X_n$
of dimension $2n-4$, obtained as a linear section of the Grassmannian $G(2,2n)$.
If $n=3$, then it is proved in \cite{Ku4} that the 2-Calabi-Yau category $\CCC_3$
is equivalent to the derived category $\D^b(X_3)$. When $n\geq 4$, $Y_n$ is singular
in codimension 5, and one has to replace $\D^b(Y_n)$ by a categorical resolution of
singularities $\tilde\D_n = \widetilde{\D^b(Y_n)}$, the noncommutative Calabi-Yau $\CCC_n$
being defined as the orthogonal complement to the exceptional collection
$\CO$, $\CO(1)$, \dots, $\CO(n-1)$ in $\tilde\D_n$. The
wanted equivalence $\CCC_n\simeq \D^b(X_n)$ remains conjectural
because of the technical difficulties
of proving the existence of $\tilde\D_n$ (Conjectures \ref{grpf}, \ref{soyv}).
The relation
between $X_n$, $\tilde\D_n$ is a part of the Homological Projective Duality
program, started in \cite{Ku2}.

6. We prove that $F(Y_n)$ is birational to a certain subvariety $H(X_n)$
of the Hilbert scheme $X_n^{[C_{n-1}]}$, where $C_{n-1}$ is a Catalan number
(Theorem \ref{FY_bir_HX}). The proof uses an explicit description
of resolutions of singularities of Pfaffian
hypersurfaces and of their Fano schemes. For $n=3$, we get $C_2=2$, and $H(X_3)=X_3^{[2]}$,
which brings us back to the result of Beauville-Donagi, and it is obvious
that the $2$-form $\alp_2$ on $F(Y_3)$ is proportional to the 2-form induced
by the holomorphic volume element of the K3 surface $X_3$.
For bigger $n$, $F(Y_n)$, $H(X_n)$ are singular and the map
between them is not biregular, so the assertion that a holomorphic volume element of $X_n$
induces $\alp_{2n-4}$ is not obvious,  and we state it as Conjecture
\ref{4-forms_compatible}.

An important new ingredient of our techniques, comparing to
the paper \cite{KM}, where the case $n=3$ was treated, is the Hochschild
(or cyclic) homology.
The Hochschild homology $\HOH(Y)$ of a smooth projective
variety $Y$ was introduced by Markarian in \cite{Ma}. The exponential
of the universal Atiyah class of $Y$ provides the HKR isomorphisms
$$
I_Y:\HOH(Y)\rar H^\fatdot(Y,\Omega_Y^\fatdot) = H^\fatdot(Y,\CC),\qquad
\HI_Y:H^\fatdot(Y,\CC) = H^\fatdot(Y,\Omega_Y^\fatdot) \rar \HOH(Y).
$$
To any pair of smooth projective varieties $Y$ and $\CM$ and an object
$\EEE$ of the derived category $\D^b(Y\times \CM)$, one can associate
the Fourier-Mukai transform (or kernel functor) $\Phi_\EEE:\D^b(Y) \to \D^b(\CM)$
and its $\HOH$-descent $\phi_\EEE:\HOH(Y) \to \HOH(\CM)$. The natural
map induced on the Dolbeault-cohomology level is $\AJ_{\ch(\EEE)}$,
where $\AJ_\xi(c) = p_*(q^*(c)\wedge\xi)$, and $p:Y\times\CM \to \CM$ and $q:Y\times\CM\to Y$ are the projections. It turns out that the maps commute with
the HKR isomorphisms only upon some modifications, which we specify
in Proposition~\ref{commutative} and Theorem\ref{msth}.
Thus the Hochschild homology of smooth projective varieties can replace
their Dolbeault cohomology.

On the other hand, Hochschild homology can be defined for any triangulated category
which is equivalent to a semiorthogonal component of the derived category
of a smooth projective variety \cite{Ku6}.
In particular, $\HOH(\CCC_n)$ is well defined for any $Y_n$ and in the Pfaffian case
$\HOH(\CCC_n)$ is expected to be isomorphic to $\HOH(X_n) = H^\fatdot(X_n,\CC)$.
Moreover, it is shown in {\em loc.\ cit.}\/ that Hochschild homology is {\em additive}\/
with respect to semiorthogonal decompositions. So, $\HOH(\CCC_n)$ is a direct
summand in $\HOH(Y_n)$, and we prove that the projection $\HOH(Y_n) \to \HOH(\CCC_n)$
takes a generator $\omega \in H^{2n-3,1}(Y_n)$ to a holomorphic volume form of $\CCC_n$.
\bigskip

As we have already mentioned, $\CCC_n$ is expected to be equivalent
to the derived category of a genuine Calabi-Yau for a Pfaffian $Y_n$.
It is interesting to find other hypersurfaces $Y_n$ with the same
property, which we call {\em derived-special}.
When $n=3$, there are non-Pfaffian derived-special cubic fourfolds
(e.g.\ one can take a fourfold containing a plane and another 2-cycle
which has odd intersection with the 2-quadric residual for the plane).
All of them are also {\em cohomologically-special} (in the sense of Hassett~\cite{Ha})
and eventually rational. It is very interesting to understand the relation between
these notions. A natural conjecture is \textcolor{black}{that} a cubic $Y_3$ is rational
if and only if it is derived-special~\cite{Ku4}.

Remark also that on one hand our construction of the $p$-form $\alp_p$ via the linkage class
is valid for any moduli space of sheaves on $Y_n$, and on the other hand, every
moduli space of sheaves on $X_n$ carries a natural $p$-form induced by the holomorphic volume element of $X_n$.
It is interesting to find other pairs of moduli spaces, related by a Fourier-Mukai transform
like $F(Y_n), H(X_n)$, and study their $p$-forms. R.~Thomas constructed
in \cite{Th} an example of a Calabi--Yau threefold $V$ and a 3-dimensional
moduli space $M$, on which the induced 3-form is nondegenerate, so that $M$
is a new Calabi--Yau threefold associated to $V$, a Fourier-Mukai partner of it.
Are there special hypersurfaces $Y_n$ and moduli spaces
on them which are (birational to) Fourier-Mukai partners of the corresponding varieties $X_n$?

We will now describe the content of the paper by sections.

In Section \ref{iajm}, we define the Abel--Jacobi map in the transcendental setting
and describe its properties for the case of the Fano scheme of lines $F(Y_n)$ on a hypersurface $Y_n$.
In particular, we use it to define the $(2n-4)$-form $\alpha_\omega$ on $F(Y_n)$.
For $n=4$, that is when $Y$ is a 6-dimensional quartic, we also determine
the relevant Hodge numbers of $Y$ and $F$ (Proposition~\ref{Hodge_num}, proved in Appendix B).
In Section \ref{ajhc} we introduce the Hochschild homology interpretation
of the Abel--Jacobi map and show that the form $\alpha_\omega$ on $F(Y_n)$
is induced by the holomorphic volume form of a noncommutative Calabi--Yau
variety associated to $Y_n$.

In Sections~\ref{phtd}--\ref{fspv} we investigate a special case of Pfaffian hypersurfaces.
In Section \ref{phtd}, we discuss the homological projective duality between the Grassmannian
of lines and the Pfaffian variety. We define the Calabi--Yau linear section $X=X_{2n-4}$
of the Grassmannian $\Gr(2,2n)$ associated to the Pfaffian hypersurface $Y=Y_n \subset \PP^{2n-1}$,
and formulate a conjectural relation between their derived categories of coherent sheaves.
In Section \ref{fspv}, it is proved that the Fano
scheme $F(Y)$ is irreducible and has a natural crepant resolution for a generic Pfaffian $Y$. A moduli space $H(X)$ of sheaves
on $X$ is constructed, which is birational to $F(Y)$.
These results base upon some general facts from geometry of lines
on Pfaffian varieties defined by rank-$2k$ Pfaffians in $\PP(\wedge^2W^*)$
($2k \leq 2n = \dim W$), and we gather such facts in Appendix~\ref{lpv}.  In particular,
we construct there a natural crepant resolution of singularities of the variety
of lines on the Pfaffian hypersurface ($k=n-1$).

In Section \ref{eflc}, we define the linkage class $\eps_\CF$, show that it factors
through the Atiyah class and prove that formula \eqref{p-form}
provides the value of $\alp_{2n-4}$. We also prove a nonvanishing result for~$\eps_\CF$.

In the last two Sections~\ref{ec} and~\ref{4-forms} we consider another special case, $n = 4$.
In Section \ref{ec}, an explicit calculation of $\alp_4$ in coordinates is done
for the quartic 6-fold $Y_4$. In the last Section \ref{4-forms}, we describe
the classification of the orbits of $GL(7)$ in $\wedge^4\CC^7$ and show that
$\alp_4$ is minimally degenerate at the generic point of the 7-fold $F(Y_4)$.

In what follows, the base field $k$ is an algebraically closed field of characteristic 0;
sometimes we assume that $k=\CC$.
A {\em variety} is a reduced irreducible (separated) scheme
of finite type over~$k$.

{\sc Acknowledgements.} The authors thank Katia Amerik, Ugo Bruzzo, Jan Nagel, and Fedor Zak for discussions, and Claire Voisin for pointing out the reference to the paper
of Shimada. D.M. acknowledges the hospitality of
the SISSA in Trieste, where he started to work on the present paper.

\section{Abel--Jacobi map of the family of lines}\label{iajm}

Let $n\geq 3$, and let $Y=Y_n$ be
a nonsingular hypersurface  of degree $n$ in $\PP^{2n-1}$.
Let $F(Y)$ be the Fano scheme of $Y$. One can think of it
either as the Hilbert scheme of lines in $Y$, or as the locus
of the Grassmannian $G=G(2,2n)$ parameterizing the lines in $\PP^{2n-1}$
that are contained in $Y$. Unlike the case of a cubic hypersurface treated by Clemens--Griffiths \cite{CG},
where the Fano scheme of lines
is smooth whenever the hypersurface is smooth,
the Fano scheme of a smooth higher degree hypersurface may be singular
and even non-reduced, see \cite{T} for the case of a quartic threefold.
Thus we will assume in the sequel
that $Y$ is generic, and this assumption is essential.

\begin{lemma}
Let $Y$ be
a generic hypersurface  of degree $n$ in $\PP^{2n-1}$. Then $F(Y)$
is a nonsingular projective variety of dimension $3n-5$.
\end{lemma}

\begin{proof}
Denote by $V$ the vector space $\CC^{2n}$ whose projectivization contains $Y$.
So $\PP^{2n-1}=\PP(V)$, and $Y=Y_f$ is defined by a homogeneous form $f\in S^dV^*$.
The Grassmannian $G=G(2,2n)$ of lines in $\PP^{2n-1}$ is naturally embedded
into $\PP(\wedge^2V)$ via the Pl\"ucker embedding and carries the tautological
subbundle $T$ such that $H^0(G,S^nT^*)$ is canonically
isomorphic to $S^nV^*$ for any $n\geq 1$.
According to \textcolor{black}{Theorem 3.3 (iv) of \cite{AK1}}, $F(Y)$ is
the scheme of zeros of the section $s_f$ of $S^dT^*$
corresponding to $f$ under the above canonical isomorphism (with $n=d$)
as soon as the latter has expected dimension.
As $T^*$ is generated by global sections, so is $S^dT^*$, and
by Kleiman's generalization of Bertini's Theorem \cite{K},
for generic $f$, the section $s_f$
is transversal to the zero section of $S^dT^*$. Thus the zero locus of $s_f$ is
smooth and is of expected codimension $n+1=\rk S^dT^*$ whenever it is nonempty.
The fact that $s_f$ has
nonempty zero locus follows from the results of Barth--Van de Ven \cite{BvdV}
and Debarre--Manivel \cite{DM}.
\end{proof}

Following \cite{BD}, we define the Abel--Jacobi map
on the middle-dimensional integer cohomology of $Y$
by the formula
\begin{equation}\label{AJmap}
\AJ : H^{2n-2}(Y,\ZZ)\lra H^{2n-4}(F,\ZZ), \ \ \ \AJ(c)=p_*q^*(c),
\end{equation}
where $F=F(Y)$, $p:Z\rar F$, $q:Z\rar Y$ are the natural projections from
the universal family of lines $Z$ in $Y$,
$$
Z=\{(\ell,x)\in F\times Y\mid x\in\ell\}.
$$

We have defined the Abel--Jacobi map on the integral cohomology.
Going over to the complex coefficients, we can represent the cohomology classes
on $Y$ by closed forms of type $(i,j)$, and $\AJ$ lifts to
the level of closed forms as
the map of integration over the fibers $p$. This allows us
to consider it as a morphism of Hodge structures, shifting the
weight by $-2$.


We will denote by $h$ the class of a hyperplane section of $Y$, and by
$\sigma_i$, $\sigma_{ij}$ the Schubert classes on
$G=\Gr(2,2n)$, as well as their restrictions to $F\subset G$.
The primitive parts of the above cohomology groups are defined
by
$$
H^{2n-2}(Y,\ZZ)_{prim}=\{x\in H^{2n-2}(Y,\ZZ)\ | \ xh=0\},$$   $$
H^{2n-4}(F,\ZZ)_{prim}=\{u\in H^{2n-4}(F,\ZZ)\ | \ u\sigma_1^{n}=0\}.
$$
We introduce also the vanishing part of the cohomology of $F$:
$$
H^{2n-4}(F,\ZZ)_{van}=\ker\big( H^{2n-4}(F,\ZZ)\xrightarrow{{\mathrm{Gysin}}}
H^{6n-4}(G,\ZZ)\big),
$$
where the Gysin map is associated to the natural embedding $F\hookrightarrow G$.
Obviously, $H^{2n-4}(F,\ZZ)_{van}\subset H^{2n-4}(F,\ZZ)_{prim}$.
We will see that though in the Beauville--Donagi case ($n=3$) this inclusion is an isomorphism, this does not hold for $n=4$, in which case $H^{4}(F,\ZZ)_{van}$
is of codimension 1 in $H^{4}(F,\ZZ)_{prim}$.
One can also define the
vanishing cohomology of $Y$ as the kernel of the Gysin map associated to the embedding
into $\PP^{2n-1}$, but it will coincide with the primitive cohomology.
The following result was proved by Shimada:

\begin{theorem}[Shimada \cite{Sh}]
The Abel--Jacobi map induces an isomorphism
$$
H^{2n-2}(Y,\ZZ)_{prim}\simeq H^{2n-4}(F,\ZZ)_{van}/{\mathrm{(torsion)}}.
$$
\end{theorem}

Now, following the approach of \cite{BD}, we will show that $AJ$ is a
morphism of {\em polarized} Hodge structures, that is, it respects
natural bilinear forms on the primitive cohomology.

Remark that $Z=\PP(T_F)$, where $T_F$ is the restriction to $F$
of the universal rank-2 bundle $T$ over $G$, and
the tautological sheaf $\OOO_{Z/F}(1)$ coincides with $q^*\OOO_{Y}(1)$.
By the theory of Chern classes, $H^\fatdot(Z,\ZZ)=p^*H^\fatdot(F,\ZZ)[q^*h]$
with a single relation
$$
(q^*h)^2=p^*\sigma_1q^*h-p^*\sigma_{11}.
$$
In particular, we can write
$$
q^*(x)=p^*x_{2n-4}q^*h-p^*x_{2n-2}
$$
for any $x\in H^{2n-2}(Y,\ZZ)$, where $x_i\in H^{i}(F,\ZZ)$ and $x_{2n-4}=\AJ(x)$.

Let us define on $H^{2n-4}(F,\ZZ)$ the symmetric bilinear form $\phi$
by
$$
\phi (u,v)=\sigma_1^{n-1}uv\ \ \forall \ \ u,\ v\in H^{2n-4}(F,\ZZ).
$$
Using the above description of $H^\fatdot(Z,\ZZ)$, one can easily prove the following
lemma (by the same arguments as in \cite{BD}):
\begin{lemma}\label{char_prim}
(i) An element $x\in H^{2n-2}(Y,\ZZ)$ is primitive if and only if
$$
x_{2n-4}\sigma_{11}=0,\ \ x_{2n-2}=x_{2n-4}\sigma_{1}.
$$
(ii) For any $x,y\in H^{2n-2}(Y,\ZZ)_{prim}$,
$$
xy=-\frac{1}{n!}\phi(\AJ(x),\AJ(y)).
$$
\end{lemma}

The intersection form $(x,y)\mapsto xy$ on $H^{2n-2}(Y,\ZZ)$ is nondegenerate,
so part (ii) gives an easy proof of the injectivity of
$\AJ$ on $H^{2n-2}(Y,\ZZ)_{prim}$. But the fact that the
image is contained in $H^{2n-4}(F,\ZZ)_{van}$ (and hence in $H^{2n-4}(F,\ZZ)_{prim}$)
and the surjectivity
modulo torsion are the nontrivial parts of Shimada's result.


\begin{proposition}\label{Hodge_isometry}
$\AJ|_{H^{2n-2}(Y,\ZZ)_{prim}}$ is a Hodge isometry onto a saturated
sublattice $H^{2n-4}(F,\ZZ)_{van}$ of $H^{2n-4}(F,\ZZ)_{prim}$,
polarized by the bilinear form $-\frac{1}{n!}\phi$.
\end{proposition}

As the image of $H^{2n-3,1}(Y)\subset {H^{2n-2}(Y,\CC)_{prim}}$
under $\AJ$ is in $H^{2n-4,0}(F)$, we obtain:

\begin{corollary}\label{2n-4}
$F$ carries a nonzero $(2n-4)$-form $\alp\in H^0(F,\Omega^{2n-4})$. It is defined uniquely
up to proportionality as a generator of the one-dimensional vector space
$\AJ(H^{2n-3,1}(Y))=H^{2n-4,0}(F)$.
\end{corollary}

\begin{proof}
The Hodge numbers of $Y$ are  given by the Griffiths residue theorem
(see e. g. \cite{G}): for a smooth degree-$d$ hypersurface $Y_f\subset \PP^N$,
the primitive Hodge cohomology $H^{N-p,p-1}(Y_d)_{\rm prim}$
is identified with the homogeneous component of degree $pd-N-1$ of the
Jacobian ring $R_f=k[x_0,\ldots,x_N]/
(\partial f/\partial x_0,\ldots,\partial f/\partial x_N)$.
Here one gets $h^{2n-3,1}(Y)=1$, and Shimada's Theorem implies
the result.\textbf{}
\end{proof}

Now we will turn to the case $n=4$. \textcolor{black}{This explicit
example is of great interest, because it is the first example which goes beyond the
Beauville-Donagi case and illustrates very well the phenomena arising in higher
dimensions: the discrepancy between the primitive and vanishing cohomology,
the degenerate $(2n-4)$-form on $F(Y_4)$,\ldots Moreover, this case is the
first to consider in the problem we suggested in the introduction, to search for
genuine Calabi-Yau's that might realize the CY category $\CCC_n$ for special $Y_n$'s.
Such a Calabi-Yau is expected to have a (piece of its) Hodge structure, isomorphic
to the Hodge structure of $F(Y_n)$, so it is worthy to compute the latter.
}

We will calculate some of the Hodge numbers of $F=F(Y_4)$, and the
results will imply that
$\rk H^{4}(F,\ZZ)_{prim}=
\rk H^{6}(Y_4,\ZZ)_{prim}+1$. The difference is due to the fact
that the ambient Grassmannian $G=\Gr(2,8)$ has a rank two
cohomology $H^{4}(G,\ZZ)$, and its restriction to $F$
is injective, so that the
Gysin map $H^{4}(F,\QQ)\rar H^8(G,\QQ)$ is surjective.

The proof of the
following proposition is given in Appendix~\ref{proof_betti_numbers}.

\begin{proposition}\label{Hodge_num}
Let $Y=Y_4$ be a generic quartic $6$-fold in $\PP^7$,
and $F=F(Y)$ its Fano scheme. Then we have the following results for
their Hodge numbers:

(i) $h^{6,0}(Y)=0,\ h^{5,1}(Y)=1,\ h^{4,2}(Y)=266, \ h^{3,3}(Y)=1108.$

(ii) $h^{0,0}(F)=1$ (that is, $F$ is connected), $h^{4,0}(F)=1$,
$h^{7,0}(F)=336$, and $h^{i,0}(F)=0$ for $i$ different from $0, 4, 7$.

(iii) $h^{1,1}(F)=1$, $h^{1,2}(F)=0$, $h^{1,3}(F)=266$, $h^{2,2}(F)=1109$.

\end{proposition}

\begin{corollary}\label{4_form}
For generic $Y=Y_4$, $F$ is a smooth connected $7$-dimensional
projective variety, and $H^0(F,\Omega^4_F)$ is $1$-dimensional,
generated by the $4$-form $\alp_\omega$, the Abel--Jacobi
image of a generator $\omega$ of $H^{5,1}(Y)\simeq \CC$.
\end{corollary}

In the sequel we will investigate this $4$-form in more detail.
In particular, we believe that $F$ is a moduli space of sheaves
on some (categorical deformation of a)
Calabi--Yau 4-fold $X$, and $\alp_\omega$ has a different
interpretation as a 4-form induced by the holomorphic volume element of $X$.
We can produce $X$ only for special quartics $Y$, namely, for Pfaffian ones.

\section{Abel--Jacobi map and Hochschild homology}
\label{ajhc}

The natural context for the Abel--Jacobi map is the Hochschild (or cyclic) homology.
The Hochschild homology of the derived category of coherent sheaves was introduced
and investigated by N.~Markaryan in his famous preprint~\cite{Ma}. Some points of his
definition was later clarified and developed by A.~C\u ald\u araru~\cite{Ca1,Ca2}, see
also~\cite{Ra} and~\cite{MS}.

We first recall the definition.

\begin{definition}[\cite{Ma}]
Let $X$ be a smooth projective variety.
The Hochschild homology of $X$ is defined as
$$
\HOH(X) = \HH^\bullet(X\times X,\Delta_*\CO_X \mathop{\otimes}\limits^\LL \Delta_*\CO_X),
$$
where $\Delta:X \to X\times X$ is the diagonal embedding,
and $\HH^\bullet$ is the hypercohomology.
\end{definition}

Note that we have the following natural identifications
\begin{equation}\label{hoh1}
\HOH(X) \cong \HH^\bullet(X,L\Delta^*\Delta_*\CO_X),
\end{equation}
and
\begin{equation}\label{hoh2}
\HOH(X) \cong \HH^\bullet(X,\Delta^!\Delta_*\omega_X[\dim X]).
\end{equation}
The first follows immediately from the definition (using the projection formula),
and the second is obtained from the first using the duality isomorphism
$$
D_X: \Delta^! \Delta_* F \cong
L\Delta^*(\Delta_*F) \otimes \omega_X^{-1}[-\dim X] \cong
L\Delta^*\Delta_*(F \otimes \omega_X^{-1}[-\dim X]).
$$

The Hochschild homology of $X$ is closely related to its Dolbeault cohomology.
To state this relation, we need the notion of the Atiyah class \cite{Ill}.

\begin{definition}
Let $X$ be an algebraic variety, $\Delta:X \to X\times X$ the diagonal embedding,
$\CI_\Delta$ the ideal sheaf
of the diagonal $\Delta(X) \subset X\times X$, and
$\Delta(X)^{(2)} \subset X\times X$ the second infinitesimal neighborhood of the diagonal, that is the closed subscheme
of $X\times X$ defined by the ideal sheaf $\CI_\Delta^2$. Then
$\CI_\Delta/\CI_\Delta^2 \cong \CN^\dual_{\Delta(X)/Y\times X} \cong \Omega^1_X$, and there is a natural
exact sequence
\begin{equation}\label{uat}
0 \to \Delta_*\Omega^1_X \to \CO_{\Delta(X)^{(2)}} \to \Delta_*\CO_X \to 0.
\end{equation}
The extension class
$$
\widetilde{\At}_X \in \Ext^1(\Delta_*\CO_X,\Delta_*\Omega^1_X).
$$
of this exact triple is called the universal Atiyah class on $X$.
Further, let $\FFF$ be a sheaf on $X$ or an object
of $\D^b(X)$, and
\begin{equation}\label{atF}
0\rar \FFF\otimes\Omega^1_X\rar \pr_{2*}(\pr_1^*\FFF\otimes \CO_{\Delta(X)^{(2)}})
\rar \FFF\rar 0
\end{equation}
the exact triple obtained by applying
$\pr_{2*}(\pr_1^*\FFF\otimes \: \fatdot\: )$ to \eqref{uat}.
The extension class $\At_\CF\in\Ext^1(\CF,\CF\otimes\Omega^1_X)$ of \eqref{atF}
is called the Atiyah class of $\FFF$, and
the object in the middle of \eqref{atF} is the sheaf (or the complex)
of first jets of $\FFF$.
\end{definition}

Consider the following maps
$$
I_X^l:
\xymatrix@1{
L\Delta^*\Delta_*\CO_X  \ar[rr]^-{\frac1{l!}\widetilde{\At_X}^{\wedge l}} &&
L\Delta^*\Delta_*\Omega^l_X[l] \ar[rr]^\lambda &&
\Omega^l_X[l],
}
$$
and
$$
\HI_X^{d-l}:
\xymatrix@1{
\Omega^l_X[l] \ar[rr]^\rho &&
\Delta^!\Delta_*\Omega^l_X[l] \ar[rrr]^-{\frac1{(d - l)!}\widetilde{\At_X}^{\wedge (d - l)}} &&&
\Delta^!\Delta_*\omega_X[\dim X],
}
$$
where $\lambda$ and $\rho$ are the natural adjunction morphisms.
Summing up these maps on $l$, we obtain the maps
$$
I_X:L\Delta^*\Delta_*\CO_X \to \oplus_k \Omega^k_X[k]
\qquad\text{and}\qquad
\HI_X: \oplus_k \Omega^k_X[k] \to \Delta^!\Delta_*\omega_X[\dim X].
$$

\begin{theorem}[\cite{Ma,Ra}]\label{hohx}
The maps
$$
I_X:\HH^\bullet(X,L\Delta^*\Delta_*\CO_X) \to \oplus_{p,q} H^q(X,\Omega^p_X),
\quad
\HI_X: \oplus_{p,q} H^q(X,\Omega^p_X) \to \HH^\bullet(X,\Delta^!\Delta_*\omega_X[\dim X])
$$
are isomorphisms. Moreover, the following diagram is commutative
$$
\xymatrix{
\oplus_{p,q} H^q(X,\Omega^p_X) \ar[rr]^-{\td_X^{-1}\circ J_X} \ar[d]_{\HI_X} && \oplus_{p,q} H^q(X,\Omega^p_X) \\
\HH^\bullet(X,\Delta^!\Delta_*\omega_X[\dim X]) \ar[rr]^-{D_X} && \HH^\bullet(X,L\Delta^*\Delta_*\CO_X) \ar[u]_{I_X}
}
$$
where $\td_X$ is the Todd genus of $X$, and $J_X$ acts on $H^q(X,\Omega^p_X)$ as the  multiplication by $(-1)^p$.
\end{theorem}

The isomorphism $I_X$ is known as the {\sf Hochschild--Kostant--Rosenberg (HKR) isomorphism}\/ and
$\hat{I}_X$ as the {\sf twisted Hochschild--Kostant--Rosenberg isomorphism}.


Another important point about Hochschild homology is its functoriality.
For any morphism of smooth algebraic varieties $f:X \to Y$, there are pullback and pushforward maps
on Hochschild homology, $f^*:\HOH(Y) \to \HOH(X)$ and $f_*:\HOH(X) \to \HOH(Y)$ respectively.
The pullback map takes any $h \in \HOH(Y) = \Hom^\bullet(\CO_Y,L\Delta^*\Delta_*\CO_Y)$ to the composition
$$
\xymatrix@1{
\CO_X \ar@{=}[r] &
Lf^*\CO_Y \ar[r]^-{h} &
Lf^*L\Delta^*\Delta_*\CO_Y \ar[r]^-{bc} &
L\Delta^*\Delta_*Lf^*\CO_Y \ar@{=}[r] &
L\Delta^*\Delta_*\CO_X
},
$$
where $bc$ stands for the base change morphism.
Similarly, the pushforward map takes any $h \in \HOH(X) = \Hom^\bullet(\CO_X,L\Delta^!\Delta_*\omega_X[d_X])$ to the composition
$$
\xymatrix@1{
\CO_Y \ar[r] &
Rf_*\CO_X \ar[r]^-{h} &
Rf_*\Delta^!\Delta_*\omega_X[d_X] \ar[r]^-{bc} &
\Delta^!\Delta_*Rf_*\omega_X[d_X] \ar[r] &
\Delta^!\Delta_*\omega_Y[d_Y]
},
$$
where $d_X=\dim X$, $d_Y=\dim Y$.
On the other hand, with any object $\CE \in \D^b(X)$ one can associate the map $\tau_\CE:\HOH(X) \to \HOH(X)$
which takes any $h \in \HOH(X) = \Hom^\bullet(\CO_X,L\Delta^*\Delta_*\CO_X)$ to the composition
%
\begin{multline*}
\CO_X \stackrel{h}\to
L\Delta^*\Delta_*\CO_X \to
\CE^\vee \mathop{\otimes}\limits^\LL \CE \mathop{\otimes}\limits^\LL L\Delta^*\Delta_*\CO_X \cong \\
L\Delta^*((\CE^\vee \boxtimes \CE) \mathop{\otimes}\limits^\LL \Delta_*\CO_X) \cong
L\Delta^*\Delta_*(\CE^\vee \mathop{\otimes}\limits^\LL \CE) \stackrel{\tr}\to
L\Delta^*\Delta_*(\CO_X)
\end{multline*}
Combining all these operations, one can associate a map $\phi_\CE:\HOH(X) \to \HOH(Y)$ of Hochschild homology
to any kernel $\CE \in \D^b(X\times Y)$ as follows:
$$
\phi_\CE = p_{Y*}\circ\tau_\CE\circ p_X^*.
$$

It turns out, however, that the functoriality of the Hochschild homology is not completely compatible
with the functoriality of the Dolbeault cohomology. Although the HKR
isomorphism commutes with pullbacks, it does not commute with pushforwards. Similarly,
though the twisted HKR isomorphism commutes with pushforwards, it does not
commute with pullbacks. To get the compatibility, we twist by the square root of the Todd genus.

Define the {\sf modified HKR isomorphism} $\MI_X:\HOH(X) \to \oplus_{p,q} H^q(X,\Omega_X)$
by
$$
\MI_X(h) = I_X(h)\wedge\td_X^{1/2}.
$$
For each $\xi \in H^\fatdot(X\times Y,\CC)$ we define the Abel--Jacobi map $\AJ_\xi$
by
\begin{equation}\label{TAJmap}
\AJ_\xi : H^\bullet(X,\CC)\lra H^{\bullet}(Y,\CC), \qquad \AJ(c)=p_{Y*}(p_X^*(c)\cdot\xi).
\end{equation}

\begin{theorem}[\cite{MS}]\label{msth}
For any $\CE \in \textcolor{black}{\D^b(X\times Y)}$, the following diagram is commutative:
$$
\xymatrix{
\HOH(X) \ar[rrrrrr]^-{\phi_\CE} \ar[d]_{\MI_X} &&&&&& \HOH(Y) \ar[d]^{\MI_Y} \\
\oplus_{p,q} H^q(X,\Omega_X)
\ar[rrrrrr]^-{{\AJ_{\ch(\CE)\wedge\td_X^{1/2}\wedge\td_Y^{1/2}}}} &&&&&&
\oplus_{p,q} H^q(Y,\Omega_Y)
}
$$
\end{theorem}

Thus, by modifying appropriately the HKR isomorphism and the Abel--Jacobi map one can get the compatibility.
However, there is another way which is more convenient for our purposes.

\begin{proposition}\label{commutative}
For any $\CE\in\D^b(X\times Y)$ there is a commutative diagram
$$
\xymatrix{
\HOH(X) \ar[rrrr]^-{\phi_{\CE}} &&&& \HOH(Y) \ar[d]^{I_Y} \\
\oplus H^p(X,\Omega^q_X) \ar[u]^{D_X\HI_X} \ar[rrrr]^{\AJ_{\ch(\CE)}\circ J_X} &&&& \oplus H^p(Y,\Omega^q_Y) \\
}
$$
\end{proposition}
\begin{proof}
Indeed,
\begin{multline*}
I_Y(\phi_\CE(D_X\HI_X(\omega))) =
\MI_Y(\phi_\CE(D_X\HI_X(\omega)))\wedge\td_Y^{-1/2} = \\
\AJ_{\ch(\CE)\wedge\td_X^{1/2}\wedge\td_Y^{1/2}}(\MI_X(D_X\HI_X(\omega)))\wedge\td_Y^{-1/2} =
\AJ_{\ch(\CE)\wedge\td_X^{1/2}}(I_X(D_X\HI_X(\omega))\wedge\td_X^{1/2}) = \\
\AJ_{\ch(\CE)}(I_X(D_X\HI_X(\omega))\wedge\td_X) =
\AJ_{\ch(\CE)}(J_X(\omega)).
\end{multline*}
Here the first equality is the definition of $\MI_Y$, the second is Theorem~\ref{msth},
the third is the projection formula plus the definition of $\MI_X$, the fourth is the projection formula,
and the last one is Theorem~\ref{hohx}.
\end{proof}

Now assume that $\CM$ is a smooth and projective component of the Hilbert scheme of
a smooth projective variety $Y$. The example of the Fano scheme we started with is
of this kind, for $F(Y)$ can be interpreted as the moduli space of sheaves
parameterizing the structure (or ideal) sheaves of lines contained in $Y$.
In this case the kernel functor $\Phi_{\CO_Z}:\D^b(Y) \to \D^b(\CM)$ given
by the structure sheaf $\CO_Z$ of the universal subscheme $Z \subset \CM\times X$,
boils down to the composition $p_*q^*$, where $p:Z \to \CM$ and $q:Z \to Y$
are the projections. So, it is natural to expect a relation of the Abel--Jacobi map given by $Z$
with the Abel--Jacobi map given by $\ch(\CO_Z)$.

\begin{lemma}\label{maj}
For any $Z \subset Y \times \CM$ and any $\omega \in H^p(Y,\Omega^q_Y)$
$$
\AJ_{\ch(\CO_Z)}(\omega) = \AJ(\omega) + \text{\sf terms of higher degree}.
$$
\end{lemma}
\begin{proof}
Indeed, $\ch(\CO_Z) = [Z] +  \text{\sf terms of higher degree}$, and $\AJ_{[Z]} = \AJ$.
\end{proof}

As a consequence we obtain the following

\begin{corollary}\label{hoch_atiyah}
Let $Y \subset \PP^{2n-1}$ be a generic hypersurface of degree $n$ and let $\textcolor{black}{\MMM=}F(Y)$ be its Fano scheme of lines.
Let $\omega$ be the generator of $H^1(Y,\Omega^{2n-3}_Y)$ and
$\alpha_\omega \in H^0(F(Y),\Omega^{2n-4}_{F(Y)})$ the corresponding holomorphic form.
Then
$$
\alpha_\omega = - I_\CM(\phi_{\CO_Z}(D_Y\HI_Y(\omega))) + \text{\sf terms of higher degree}.
$$
\end{corollary}

This Hochschild homology interpretation allows us to show that
the form $\alpha_\omega$ actually comes from a certain Hochschild homology class
in a certain triangulated subcategory of $\D^b(Y)$. For this we have to discuss
the Hochschild homology of triangulated categories.

It is well known (see e.g.~\cite{Or}, Theorem 2.1.8) that $\HOH(X) \cong \HOH(Y)$
if the derived categories of $X$ and $Y$ are equivalent. Thus Hochschild
homology is an invariant of the derived category and it is natural to generalize
the definition to any triangulated category. Actually, we will need this in a particular case.

Let $\CA \subset \D^b(X)$ be an admissible subcategory (i.e. a component of a semiorthogonal
decomposition of $\D^b(X)$). It is proved in~\cite{Ku5} that there exists an object
$K(\CA) \in \D^b(X\times X)$ such that the corresponding kernel functor
$\Phi_{K(\CA)}:\D^b(X) \to \D^b(X)$ is the projection onto the subcategory $\CA$.

\begin{definition}[\cite{Ku6}]
Let $\CA \subset \D^b(X)$ be an admissible subcategory.
The Hochschild homology of $\CA$ is defined as
$$
\HOH(\CA) = \HH^\bullet(X\times X,K(\CA) \mathop{\otimes}\limits^\LL K(\CA)^*),
$$
where $K(\CA)^*$ is the kernel of the left adjoint functor to the projection onto the subcategory $\CA$.
\end{definition}

A very important property of Hochschild homology is the following additivity Theorem.

\begin{theorem}[\cite{Ku6}]\label{hohsum}
If $\D^b(X) = \langle \CA_1,\CA_2,\dots,\CA_n \rangle$ is a semiorthogonal decomposition, then
$$
\HOH(X) = \HOH(\CA_1)\oplus\HOH(\CA_2) \oplus \dots \oplus \HOH(\CA_n).
$$
Moreover, if $K \subset \CA_i\boxtimes\D^b(Y) \subset \D^b(X\times Y)$, then the map $\phi_K:\HOH(X) \to \HOH(Y)$
factors through the projection to the summand $\HOH(\CA_i)$.
\end{theorem}

Recall that a triangulated category $\CT$ is called an $m$-Calabi--Yau category,
if the shift by $m$ functor $[m]$ is a Serre functor in $\CT$, i.e.\ if we are
given bifunctorial isomorphisms
$$
\Hom(F,G) \cong \Hom(G,F[m])^*
$$
for all $F,G \in \CT$.

\begin{remark}
One can also define the Hochschild cohomology $\HOH^\bullet(\CT)$ of an admissible subcategory $\CT \subset \D^b(X)$.
The Hochschild cohomology is naturally a ring, while the Hochschild homology is a module over it.
\end{remark}

A triangulated category $\CT$ is called {\em connected}\/ if $\HOH^{<0}(\CT) = 0$ and $\HOH^0(\CT) = \CC$.
One has $\HOH^0(\D^b(X)) = H^0(X,\CO_X)$, so the derived category of a smooth projective variety
is connected if and only if the variety itself is connected.

\begin{lemma}[\cite{Ku6}]
Let $\CT$ be a connected $m$-Calabi--Yau category, equivalent to an admissible subcategory in the derived category of a smooth projective variety,
then $\HOH_i(\CT) = 0$ unless $-m \le i \le m$.
Moreover $\dim\HOH_{-m}(\CT) = 1$.
\end{lemma}

Let $\CT$ be a connected $m$-Calabi--Yau category and let $\alpha = \alpha_\CT$ be the generator of the space $\HOH_{-m}(\CT)$.
By the above Lemma, $\alpha_\CT$ is defined uniquely up to a constant. We will call $\alpha_\CT$ the {\sf holomorphic volume form of $\CT$}.

\begin{remark}
If $\CT$ is an $m$-Calabi--Yau category, then the Hochschild homology is a free module over the Hochschild cohomology ring,
and $\alpha_\CT$ is a generator.
\end{remark}

Now we are coming back to our original setting: $Y=Y_n$ is
a general nonsingular hypersurface  of degree $n$ in $\PP^{2n-1}$,
$F=F(Y)$ is its Fano scheme, $q:Z\rar Y$ is
the universal family of lines in $Y$, $p:Z\rar F$
the natural projection.

\begin{theorem}[\cite{Ku1}]\label{cccn}
The line bundles $\CO_Y,\CO_Y(1),\dots,\CO_Y(n-1)$ form an exceptional collection in $\D^b(Y)$,
so that we have a semiorthogonal decomposition
$$
\D^b(Y_n) = \langle \CCC_n,\CO_Y,\CO_Y(1),\dots,\CO_Y(n-1) \rangle,
$$
where $\CCC_n = \{ F \in \D^b(Y)\ |\ \HH^\bullet(Y,F) = \HH^\bullet(Y,F(-1)) = \dots = \HH^\bullet(Y,F(1-n)) = 0\}$.
Moreover, $\CCC_n$ is a connected $(2n-4)$-Calabi--Yau category.
\end{theorem}


\begin{corollary}
We have
$$
\HOH_i(\CCC_n) = \begin{cases}
H^{n-1+t,n-1-t}(Y,\CC), & \text{if $i = 2t \ne 0$,}\\
H^{n-1,n-1}(Y,\CC) \oplus \CC^{n-2}, & \text{if $i = 0$,}\\
0, & \text{otherwise.}
\end{cases}
$$
\end{corollary}
\begin{proof}
The category generated by an exceptional bundle is equivalent to the derived category of vector spaces,
that is, to the derived category of coherent sheaves on a point. Thus, $\HOH_\bullet(\langle \CO_Y(k) \rangle) = \CC$ (sitting in degree $0$)
for any $k$. Combining with Theorems~\ref{hohsum} and~\ref{cccn}
we conclude that
$$
\HOH_i(\D^b(Y)) = \begin{cases}
\HOH_i(\CCC_n) \oplus \CC^n, & \text{if $i = 0$}\\
\HOH_i(\CCC_n), & \text{otherwise.}
\end{cases}
$$
On the other hand, by Theorem~\ref{hohx} we have
$$
\HOH_i(\D^b(Y)) = \begin{cases}
H^{n-1+t,n-1-t}(Y,\CC), & \text{if $i = 2t \ne 0$,}\\
H^{n-1,n-1}(Y,\CC) \oplus \CC^{2n-2}, & \text{if $i = 0$,}\\
0, & \text{otherwise.}
\end{cases}
$$
So, the Corollary follows.
\end{proof}

Our final goal in this section is
to show that the form $\alpha_\omega$ on $F(Y)$
constructed above is induced by the form $\alpha_{\CCC_n} \in \HOH_{4-2n}(\CCC_n)$.

First of all choose a line  $\ell = \PP(T) \subset \PP(V)$ on $Y \subset \PP(V)$.
Then $\ell$ is the zero locus of a regular section of a vector bundle $V/T\otimes\CO_{\PP(V)}(1)$
on $\PP(V)$, hence on $\PP(V)$ we have the following Koszul resolution
$$
\left\{\Lambda^{2n-2}T^\perp\otimes\CO_{\PP(V)}(2-2n) \to \dots \to
\Lambda^{2}T^\perp\otimes\CO_{\PP(V)}(-2) \to
T^\perp\otimes\CO_{\PP(V)}(-1) \to \CO_{\PP(V)}\right\}  \cong \CO_\ell.
$$
Restricting the resolution to $Y$, twisting by $\CO_Y(k-2)$ and truncating at the term $\Lambda^{k-2}T^\perp\otimes\CO_Y$,
we obtain the complex
\begin{equation}\label{trunc_koszul}
\left\{\Lambda^{k-2}T^\perp\otimes\CO_Y \to \dots \to
T^\perp\otimes\CO_Y(k-3) \to \CO_Y(k-2) \to \CO_\ell(k-2) \right\}
\end{equation}
which we consider as an object of $\D^b(Y)$ and denote by $\CG^k_\ell$.

\begin{remark}
Note that the restriction of the Koszul complex to $Y$ has two cohomology sheaves,
$\CO_\ell$ at the rightmost term and $L_1i^*i_*\CO_\ell \cong \CO_\ell(-n)$ at the term  next to it.
Therefore, $\CG^k_\ell$ is indeed a complex with two cohomology sheaves,
a torsion free sheaf
at the leftmost term and $\CO_\ell(k-n-2)$ at the term third from the right.
\end{remark}

\begin{proposition}\label{cgcc}
We have $\CG^k_\ell \in \langle \CO_Y,\CO_Y(1), \dots, \CO_Y(k-1) \rangle^\perp$ for $k \le n$.
In particular, $\CG^n_\ell \in \CCC_n$.
\end{proposition}
\begin{proof}
First of all, note that $\CG^1_\ell = \CO_\ell(-1)$ is right orthogonal to $\CO_Y$, since we have
$\Ext^\fatdot(\CO_Y(1),\CO_\ell) = H^\fatdot(Y,\CO_\ell(-1)) = 0$.
On the other hand, we have an exact triangle
$$
\CG^k_\ell \to \Lambda^{k-2}T^\perp\otimes\CO_Y \to \CG^{k-1}_\ell(1)
$$
It follows by induction that $\CG^k_\ell \in \langle \CO_Y(1), \dots, \CO_Y(k-1) \rangle^\perp$,
so it remains to check that $\Ext^\fatdot(\CO_Y,\CG^k_\ell) = H^\fatdot(Y,\CG^k_\ell) = 0$.
For this we apply the functor $H^\fatdot(Y,-)$ to the complex~\eqref{trunc_koszul}.
We get a complex
$$
\Lambda^{k-2}T^\perp \to \Lambda^{k-3}T^\perp \otimes V^* \to \dots \to
\Lambda^{2}T^\perp\otimes S^{k-4}V^* \to
T^\perp\otimes S^{k-3}V^*  \to S^{k-2}V^*  \to S^{k-2}T^*
$$
which is well known to be exact. So, $H^\fatdot(Y,\CG^k_\ell) = 0$ and we are done.
\end{proof}

Denote by $\LL$ the endofunctor of the category $\D^b(Y)$ defined as the composition
of a twist by $\CO_Y(1)$ and a left mutation through $\CO_Y$. In other words,
for any $F \in \D^b(Y)$ there is an exact triangle
$$
H^\fatdot(Y,F\textcolor{black}{(1)})\otimes\CO_Y \to F(1) \to \LL(F).
$$
Then $\CG^k_\ell = \LL^{k-1}(\CO_\ell(-1))$. Moreover, one can check (see~\cite{Ku1}) that $\LL$ induces an autoequivalence
of the category $\CCC_n$ such that $\LL^n \cong [2 - n]$ on $\CCC_n$. On the other hand, repeating the arguments
of Proposition~5.4 of~\cite{KM} one can check that $\LL$ induces a chain of isomorphisms
\begin{equation}\label{many_iso}
\Ext^1(\CO_\ell,\CO_\ell) =
\Ext^1(\CO_\ell(-1),\CO_\ell(-1)) \to
\Ext^1(\CG^2_\ell,\CG^2_\ell) \to \dots \to
\Ext^1(\CG^n_\ell,\CG^n_\ell).
\end{equation}
It follows that the form defined on $F(Y)$ by the structure sheaf $\CO_Z$ of the universal line coincides
(up to a sign) with the form defined by the family $\{\CG^n_\ell\}_{\ell \in F(Y)}$.
\textcolor{black}{More
precisely, we can replace, $\phi_{\OOO_Z}$ by $\pm\phi_{\CG^n}$, where $\CG^n = ({\mathsf{id}_{D^b(F(Y))}}\times{\mathbb{L}})^{n-1}(\CO_Z(0,-1)) \in D^b(F(Y)\times Y)$, the
functor ${\mathsf{id}_{D^b(F(Y))}}\times{\mathbb{L}}$ being given by the kernel
$\Delta_{F(Y)*}\CO_{F(Y)} \boxtimes K_{\mathbb{L}} \in
D^b(F(Y)\times F(Y)\times Y\times Y)$, where $K_{\mathbb{L}}$ denotes
the kernel providing the functor ${\mathbb{L}}$.}

On the other hand, since $\CG^n_\ell \in \CCC_n$ for any $\ell$, it follows that
$\phi_{\CG_n}:\HOH(Y) \to \HOH(F(Y))$ factors through the projection $\HOH(Y) \to \HOH(\CCC_n)$.

\begin{proposition}\label{induced_cy}
The form $\alpha_\omega \in H^0(F(Y),\Omega^{2n-4}_{F(Y)})$ is induced by the holomorphic volume form $\alpha_{\CCC_n}$
of the Calabi--Yau category $\CCC_n$.
\end{proposition}
\begin{proof}
Recall that $\alpha_\omega = \AJ(\omega)$, where $\omega$ is the generator of $H^1(Y,\Omega^{2n-3}_Y)$.
Let us show that
\begin{equation}\label{alphacg}
\alpha_\omega = \AJ_{\ch(\CG^n)}(\omega) + \text{\sf terms of higher degree}.
\end{equation}
Indeed, by definition of $\CG^n$ we have
$$
\ch(\CG^n) = \ch(\CO_Z((n-2)H)) - \sum_{p=0}^{n-2} (-1)^p\ch(\Lambda^pT^\perp)\ch(\CO_Y((n-2-p)H)),
$$
where $H$ is the class of a hyperplane section of $Y \subset \PP^{2n-1}$.
Note that $p_*q^*(\omega \wedge [H]^k) = 0$ for any $k$, hence the second summand does not
contribute into $\At_{\ch(\CG^n)}(\omega) = p_*(q^*(\omega\wedge\ch(\CG^n)))$, so
$\At_{\ch(\CG^n)}(\omega) = \At_{\ch(\CO_Z((n-2)H))}(\omega)$.
But $\ch(\CO_Z((n-2)H)) = [Z] + \text{\sf terms of higher degree}$, whence the claim.
Further, note that by Proposition~\ref{commutative} it follows that
$$
\alpha_\omega = - I_{F(Y)}(\phi_{\CG_n}(D_Y\HI_Y(\omega))) + \text{\sf terms of higher degree}.
$$
But $\CG^n \subset \CCC_n \boxtimes \D^b(F(Y))$ by Proposition~\ref{cgcc}. So, by
Theorem~\ref{hohsum} the map $\phi_{\CG_n}$ factors as $\phi'\circ\pi$, where $\pi:\HOH(Y) \to \HOH(\CCC_n)$
is the projection and $\phi':\HOH(\CCC_n) \to \HOH(F(Y))$. Since we have $\pi(D_Y\HI_Y(\omega)) \in \HOH_{4-2n}(\CCC_n) = \CC \alpha_{\CCC_n}$,
we conclude that for some $\lambda \in \CC$ we have $\alpha_\omega = \lambda I_{F(Y)}(\phi'(\alpha_{\CCC_n}) +  \text{\sf terms of higher degree}$.
\end{proof}

\textcolor{black}{
\begin{remark}\label{moduli_complexes}
A general definition of moduli spaces of complexes in a derived category is still missing
today, essentially because no general way to construct them is known, like the GIT-quotient construction. But $F(Y_n)$ has several features which enable us to consider it as a fine moduli space of objects of $\CCC_n$: (a) it parametrizes the isomorphism classes
of the complexes $\CG^n_\ell\in \mathrm{Ob}(\CCC_n)\ ({\ell \in F(Y)})$;
(b) there is a universal family $\CG^n$ over $Y_n\times F(Y_n)$, which can be
constructed by a universal version of the complex~\eqref{trunc_koszul};
(c) as follows from \eqref{many_iso}, this universal family is versal at any
point $\ell\in F(Y_n)$, that is, its Kodaira-Spencer map is an isomorphism
$T_\ell F(Y_n)\rar \Ext^1(\CG^n_\ell,\CG^n_\ell)$; (d) the same arguments, as those
proving \eqref{many_iso}, show that $\Ext^2(\CG^n_\ell,\CG^n_\ell)=0$, so that
local analytical moduli spaces of the objects $\CG^n_\ell$ exist as germs of
$\Ext^1(\CG^n_\ell,\CG^n_\ell)$ at $0$.
These local analytical moduli spaces fit into one projective variety $F(Y_n)$,
which is thus nothing else but a global moduli space.
\end{remark}
}

\section{Pfaffian hypersurfaces and their duals}\label{phtd}

Now we consider a very special type of degree $n$ hypersurfaces $Y_n \subset \PP^{2n-1}$.
A hypersurface $Y_n$ is called {\em Pfaffian}\/ if its equation can be written as the Pfaffian
of a $2n$-by-$2n$ matrix of linear forms. Such hypersurfaces are not generic, in particular
they are singular for $n \ge 4$. However, they deserve special consideration because
an analogue of the Calabi--Yau category $\CCC_n$ for these hypersurfaces has a geometric
interpretation. Actually, it is equivalent to the derived category of a certain Calabi--Yau
linear section of the Grassmannian $\Gr(2,2n)$. The relation between Pfaffian hypersurfaces
and linear sections of the Grassmannian comes naturally in the context of homological projective duality~\cite{Ku2}.

Homological projective duality (HP-duality for short) is a certain duality on the set
of smooth (non-commutative) varieties equipped with a map into a projective space
and a compatible semiorthogonal decomposition of its derived category
(called Lefschetz decomposition). It associates to a smooth variety $X$ with a map $f:X \to \PP^N$
and a Lefschetz decomposition $\CA_\fatdot$, a smooth variety $Y$ with a map into the dual projective
space $g:Y \to \check{\PP}{}^N$ and the dual Lefschetz decomposition $\CB_\fatdot$.
Classical projective duality can be considered as a quasiclassical limit of
HP-duality
since in the case where the map $f:X \to \PP^N$ is a closed embedding, the classically projectively dual
variety $X^\vee \subset \check{\PP}{}^N$ coincides with the set of critical values of the map $g:Y \to \check{\PP}{}^N$
from the HP-dual variety $Y$. But the most important property of  HP-duality is a strong relation
between the derived categories of linear sections of the HP-dual varieties $X$ and $Y$.
Choose a projective subspace $\PP(V) \subset \check{\PP}{}^N$ and let $\PP(V^\perp) \subset \PP^N$
be its orthogonal complement. Let $X_V = X\times_{\PP^N}\PP(V^\perp)$ and $Y_V = Y\times_{\check{\PP}{}^N}\PP(V)$
be the corresponding linear sections. Then the derived categories $\D^b(X_V)$ and $\D^b(Y_V)$ have
semiorthogonal decompositions, one part of each coming from the Lefschetz decomposition $\CA_\fatdot$
(or $\CB_\fatdot$) of the ambient variety (this part is called {\em trivial}), and the other parts ({\em nontrivial})
being equivalent. In some cases one of the trivial parts is zero, so the nontrivial part coincides
with the whole derived category and one obtains a fully faithful embedding of one of categories
$\D^b(X_V)$ or $\D^b(Y_V)$ into the other, with the orthogonal complement being trivial (in the above sense).

Though on the categorical level one can always describe the HP-dual variety $Y$ for any $X$,
it is a very difficult problem to find a geometrical description for it. There are not so many
examples for which the answer is known. However, Pfaffian varieties are among
them.

Let $W$ be a vector space of dimension $2n$. Consider the space $\PP(\Lambda^2W^*)$
of skew-forms on~$W$. Let $\Pf(W^*) \subset \PP(\Lambda^2W^*)$ denote the hypersurface
of degenerate skew-forms. It is called {\em the Pfaffian hypersurface}.
More generally, for each $1 \le k\le n-1$ let $\Pf_k(W^*)$ denote the subvariety
of~$\PP(\Lambda^2W^*)$ consisting of skew-forms of rank less or equal to $2n-2k$.
These varieties are called the {\em generalized Pfaffian varieties}.
They form a chain
$$
\Gr(2,W^*) = \Pf_{n-1}(W^*) \subset \Pf_{n-2}(W^*) \subset \dots \subset \Pf_2(W^*) \subset \Pf_1(W^*) = \Pf(W^*).
$$
It is easy to see that for $k < n-1$ the Pfaffian variety $\Pf_k(W^*)$ is singular,
its singular locus being the next Pfaffian variety $\Pf_{k+1}(W^*)$.
Moreover, it is easy to see that $\Pf_k(W^*)$ is the closure of a~$\GL(W)$-orbit
on $\PP(\Lambda^2W^*)$, and that $\Pf_k(W^*)$ is the $(n-2-k)$-th secant variety of $\Gr(2,W^*)$.

The sheaf of ideals of $\Pf_k(W^*)$ is the image of the map
\begin{equation}\label{sik}
\xymatrix@1{\Lambda^{2(k-1)}W \otimes \CO_{\PP(\Lambda^2W^*)}(-(n-k+1)) \ar[r]^-{\sigma_k} & \CO_{\PP(\Lambda^2W^*)}}
\end{equation}
given by Pfaffians of principal $2(n-k+1)\times2(n-k+1)$-minors
of a skew-form. Alternatively, $\sigma_k$ can be described as the unique
$\GL(W)$-semiinvariant element in the space
$$
H^0(\PP(\Lambda^2W^*),\Lambda^{2(k-1)}W^* \otimes \CO_{\PP(\Lambda^2W^*)}(n-k+1)) =
\Lambda^{2(k-1)}W^* \otimes S^{n-k+1}(\Lambda^2W^*)
$$

Another interesting fact is that the class of Pfaffian varieties is self dual with respect
to projective duality. More precisely, it is easy to see that $(\Pf_k(W^*))^\vee = \Pf_{n-1-k}(W)$.
This classical statement has the following extension.

\begin{conjecture}[\cite{Ku4}]\label{grpf}
The Pfaffian varieties $\Pf_k(W^*)$ admit categorical resolutions of singularities $\widetilde{\Pf}_k(W^*)$
such that $\widetilde{\Pf}_k(W^*)$ is Homologically Projectively Dual to $\widetilde{\Pf}_{n-1-k}(W)$.
In particular, the Grassmannian of lines $\Gr(2,W)$ is Homologically Projectively Dual to
a certain categorical resolution of singularities $\widetilde{\Pf}(W^*)$ of the Pfaffian hypersurface $\Pf(W^*)$.
\end{conjecture}

This conjecture was proved in \cite{Ku4} for $n = 3$.


Now consider the derived categories of linear sections of $X = \Gr(2,W)$ and of $\TY = \widetilde{\Pf}(W^*)$
corresponding to a generic linear subspace $V \subset \Lambda^2W^*$ of dimension $\dim V = 2n$.
Then $X_V$ is a linear section of the Grassmannian $X$ and it is easy to see that its canonical
class is zero. On the other hand, $\TY_V$ is a categorical resolution of a degree $n$ hypersurface
$Y_V \subset \PP(V) = \PP^{2n-1}$. If $\PP(V)$ does not intersect the singular locus $\Sing\Pf(W^*)$
(which is possible, by dimension reasons, only for $n \le 3$),
then $\TY_V = Y_V$, and in other cases
the resolution $\TY_V$ of $Y_V$ is nontrivial. Both $X_V$ and $\TY_V$ being linear sections of HP-dual varieties
come with semiorthogonal decompositions of their derived categories.
The one for $\D^b(X_V)$ turns out to be very simple: in this case there is no trivial part,
so the nontrivial part coincides with the whole category $\D^b(X_V)$.
As for $\D^b(Y_V)$, the trivial part is present here and is given by the exceptional
collection $\CO, \CO(1), \dots, \CO(n-1)$. So, in using the properties of HP-duality
one can deduce from Conjecture~\ref{grpf} the following
\begin{conjecture}\label{soyv}
Let $V \subset \Lambda^2W^*$ be a vector subspace of dimension $\dim V = 2n = \dim W$.
Let $V^\perp \subset \Lambda^2W$ be the orthogonal. Then the derived category
of coherent sheaves on a categorical resolution of singularities $\TY_V$ of the Pfaffian hypersurface
$Y_V = \PP(V) \cap \Pf(W^*) \subset \PP(V)$ has a semiorthogonal decomposition
with nontrivial component equivalent to the derived category of the dual linear
section of the Grassmannian $X_V = \PP(V^\perp) \cap \Gr(2,W)$. More precisely,
$$
\D^b(\widetilde{Y_V}) = \langle \D^b(X_V),\CO,\CO(1),\dots,\CO(n-1) \rangle.
$$
\end{conjecture}

See \cite{Ku4} in case $n = 3$.

If the above Conjectures are true then for a Pfaffian hypersurface $Y_n$
the Calabi--Yau category
$\CCC_n = \langle \CO,\CO(1),\dots,\CO(n-1) \rangle^\perp \subset \D^b(\TY_n)$
is equivalent to $\D^b(X_n)$. Thus, for any smooth $Y_n$ the category $\CCC_n$
can be considered as a {\em noncommutative deformation} of the derived category of the Calabi--Yau variety $X_n$.
In case $n = 3$ there are examples of other special cubics $Y_3$ for which $\CCC_3$ is
equivalent to the derived category of a (commutative) K3-surface.
It is a fascinating problem to find other special $Y_n$'s
(preferably smooth ones), for which $\CCC_n$ becomes the derived category
of a (commutative) Calabi--Yau manifold.

The semiorthogonal decomposition of Conjecture~\ref{soyv} suggests that any moduli space
of sheaves on~$Y_V$ can be represented as a moduli space on $X_V$. In the next two
sections, we will show that this
is the case for the Fano scheme of lines on $Y_V$.

\section{Fano scheme of a Pfaffian variety}
\label{fspv}


Let us fix the notation for this section.
We let $W = \CC^{2n}$, $V = \CC^{2n}$, consider a linear embedding $V \to \Lambda^2W^*$,
and denote by $X = X_V = \Gr(2,W) \cap \PP(V^\perp)$ and $Y = Y_V = \PP(V) \cap \Pf(W^*)$
the corresponding linear sections.

The goal of this section is to show that the Fano scheme
of lines on $Y = Y_V$ can be interpreted as a certain moduli space on $X = X_V$.
We start by considering the Fano scheme $F(\Pf(W^*))$ of lines on $\Pf(W^*)$.
As it is proved in Appendix~\ref{lpv} the variety
$$
\TF_1 = \{(U,L) \in \Gr(n+1,W)\times\Gr(2,\Lambda^2W^*)\ |\ L \subset \Ker(\Lambda^2W^* \to \Lambda^2U^*) \}.
$$
is a resolution of~$F(\Pf(W^*))$.

Let $F(Y)$ be the Fano scheme of lines on $Y$ and put
$\TF(Y) = F(Y)\times_{F(\Pf(W^*))}\TF_1$, which we call the resolved Fano scheme of lines on $Y$.
Let $F_0(Y) \subset F(Y)$ be the open subscheme
consisting of lines which do not intersect $Y \cap \Pf_2(W^*) \subset Y$.
Then the projection $\pi:\TF_1 \to F(\Pf(W^*))$ restricts to a projection
$\pi:\TF(Y) \to F(Y)$ which is an isomorphism over $F_0(Y)$ by the proof of
Proposition~\ref{pibir}.

\begin{lemma}
Let $n \ge 3$.
For generic $V \subset \Lambda^2W^*$ the resolved Fano scheme $\TF(Y)$ of lines on $Y = \PP(V)\cap\Pf(W^*)$
is smooth and connected and $F_0(Y)$ is nonempty. In particular, $F(Y)$ is irreducible,
$F_0(Y)$ is dense in $F(Y)$ and $\pi:\TF(Y) \to F(Y)$ is birational.
\end{lemma}
\begin{proof}
To check the smoothness of $F(Y)$ we consider the universal Pfaffian variety
and its resolved Fano scheme of lines. In other words, we consider $\Gr(2n,\Lambda^2W^*)$ and
$$
\CY = \PP_{\Gr(2n,\Lambda^2W^*)}(\CV) \times_{\PP(\Lambda^2W^*)} \Pf(W),
\quad\text{\quad}
\tilde\CF(\CY) = \Gr_{\Gr(2n,\Lambda^2W^*)}(2,\CV) \times_{\Gr(2,\Lambda^2W^*)} \TF(\Pf(W)),
$$
where we denote the tautological bundle on $\Gr(2n,\Lambda^2W^*)$ by $\CV$.
We have canonical projections $\CY \to \Gr(2n,\Lambda^2W^*)$ and $\tilde\CF(\CY) \to \Gr(2n,\Lambda^2W^*)$
and it is clear that their fibers over $V \in \Gr(2n,\Lambda^2W^*)$ are the corresponding
Pfaffian variety and its resolved Fano scheme of lines.

Note that $\tilde\CF(\CY)$ is smooth. Indeed, considering the projection $\tilde\CF(\CY) \to \TF(\Pf(W^*))$
we see that its fiber over a point $(U,L) \in \TF(\Pf(W^*))$ is just the set of all $V \in \Gr(2n,\Lambda^2W^*)$
such that $L \subset V \subset \Lambda^2W^*$. In other words, $\tilde\CF(\CY) = \Gr_{\TF(\Pf(W^*))}(2n-2,\Lambda^2W^*/L)$.

The smoothness of $\TF(Y)$ and nonemptiness of $F_0(Y)$ for general $Y$ follow
because the general fiber of a morphism of smooth varieties is smooth
and has a nontrivial intersection with a given open subset.

Now let us verify the connectedness of $\TF(Y)$.
Since the fibers of the projection $\TF(Y) \to F(Y)$ are connected,
it suffices to check that $F(Y)$ is connected. And for this it suffices
to check that $H^0(F(Y),\CO_{F(Y)}) = \CC$. But $F(Y) \subset \Gr(2,V)$
is the zero locus of a regular section of the vector bundle $S^nL^*$,
so we have the following Koszul resolution
$$
\dots \to \Lambda^2(S^nL) \to S^nL \to \CO \to \CO_{F(Y)} \to 0.
$$
Looking at the hypercohomology spectral sequence it is easy to note that
it suffices to check that $H^q(\Gr(2,V),\Lambda^t(S^nL)) = 0$ for $q \le t$ and $t > 0$.
But the Bott--Borel--Weil theorem implies that $H^q(\Gr(2,V),\Lambda^t(S^nL))$ can be nontrivial
only for $q = \dim V - 2 = 2n - 2$ and $q = 2(\dim V - 2) = 4n - 4$.
On the other hand $\rk(S^nL) = n + 1$, so $t \le n + 1$.
Since for $n \ge 3$ we have $2n - 2 \ge n + 1$ we see that
$H^q(\Gr(2,V),\Lambda^t(S^nL)) = 0$ for $q \le t$ and $t > 0$ unless
$n = 3$ and $q = t = 4$. In the latter case
$H^4(\Gr(2,6),\Lambda^4(S^3L)) = H^4(\Gr(2,6),\Sigma^{6,6}L) = H^4(\Gr(2,6),\CO(-6)) = 0$
(in this case the nontrivial cohomology is $H^8 = \CC \ne 0$), so we conclude that $H^0(F(Y),\CO_{F(Y)}) = \CC$
and that $\TF(Y)$ is connected.

Since $\TF(Y)$ is smooth and connected we conclude that $\TF(Y)$ is irreducible.
Therefore the Fano scheme $F(Y) = \pi(\TF(Y))$ is irreducible as well. And since $F_0(Y) \ne \emptyset$
the map $\pi:\TF(Y) \to F(Y)$ is birational.
\end{proof}

Denote by $\rho$ the projection $\TF(Y) \to \Gr(n+1,W)$. Our further goal is to identify
the image of~$\TF(Y)$ in $\Gr(n+1,W)$.
Consider the map $\varphi:V\otimes\CO_{\Gr(n+1,W)} \to \Lambda^2W^*\otimes\CO_{\Gr(n+1,W)} \to \Lambda^2U^*$
on the Grassmannian $\Gr(n+1,W)$. Denote
$$
\begin{array}{lll}
G(Y) &=& \{ U \in \Gr(n+1,W)\ |\ \rank(\varphi_U) \le 2n-2 \},\\
G_0(Y) &=& \{ U \in \Gr(n+1,W)\ |\ \rank(\varphi_U) = 2n-2 \}.
\end{array}
$$

\begin{lemma}
We have $\rho(F_0(Y)) = G(Y)$. Moreover, the projection
$\rho:\TF(Y) \to G(Y)$ is an isomorphism over $G_0(Y)$.
\end{lemma}
\begin{proof}
For any point $(U,L) \in \TF(Y)$ we have $\varphi_{g(L)}(L) = 0$ by definition of $\TF$.
Hence at any such point the rank of $\varphi$ is less than or equal to $\dim V/L = 2n - 2$.
On the other hand, assume that $\rank\varphi \le 2n-2$ at a point $U \in \Gr(n+1,W)$. Then there is
a two-dimensional subspace $L \subset V$ such that $\varphi_U(L) = 0$, which means that
$(U,L) \in \TF(Y)$.
\end{proof}

We have the following diagram
$$
\xymatrix{
& \TF(Y) \ar[dl]_\pi \ar[dr]^\rho \\
F_0(Y) \subset F(Y) \qquad && \qquad G(Y) \supset G_0(Y)
}
$$
where both maps $\pi$ and $\rho$ are isomorphisms over open subsets $F_0(Y)$ and $G_0(Y)$ respectively.

\begin{lemma}
For generic $Y$ the set $G_0(Y)$ is nonempty.
In particular, for generic $Y$ the projection $\rho:\TF(Y) \to G(Y)$ is birational.
\end{lemma}
\begin{proof}
Consider the universal versions of $G(Y)$ and $G_0(Y)$:
$$
\begin{array}{lll}
\CG(\CY) &=& \{ (V,U) \in \Gr(2n,\Lambda^2W^*)\times\Gr(n+1,W)\ |\ \rank(\varphi:V \to \Lambda^2U^*) \le 2n - 2 \},\\
\CG_0(\CY) &=& \{ (V,U) \in \Gr(2n,\Lambda^2W^*)\times\Gr(n+1,W)\ |\ \rank(\varphi:V \to \Lambda^2U^*) = 2n - 2 \}.
\end{array}
$$
It is easy to see that $\CG(Y)$ is irreducible and $\CG_0(Y)$ is open in $\CG(Y)$.
Hence the general fiber of $\CG_0(Y)$ over $\Gr(2n,\Lambda^2W^*)$, which is nothing
but $G_0(Y)$, is nonempty.
\end{proof}

Our further goal is to identify $G(Y)$ for $Y = Y_V$ in terms of the corresponding Calabi--Yau linear section $X = X_V$
of the Grassmannian $\Gr(2,W)$. Let
$$
d_n = \deg \Gr(2,n+1),
\qquad
c_n = (n^2 - 3n + 4)/2 = (n-1)(n-2)/2 +1
$$
Actually, $d_n$ is the Catalan number $C_{n-1}$, but we do not need it.

Consider the Hilbert scheme $\Hilb_{d_n}(X)$ of Artinian subschemes of $X$ of length $d_n$
and its subvariety
$$
H(X) = \left\{Z\in\Hilb_{d_n}(X)\ \left|\
\begin{array}{l}
\rank(W^* \to H^0(Z,\CS^*_{|Z})) \le n+1,\\
\rank(\Lambda^2W^* \to H^0(Z,\Lambda^2\CS^*_{|Z})) \le c_n
\end{array}
\right.\right\}.
$$
where $\CS$ is the tautological bundle on $G(2,W)$. If for a point
$Z \in H(X)$ the rank of the map $W^* \to H^0(Z,\CS^*_{|Z})$ equals $n + 1$,
we can associate to $Z$ the kernel of the map $W^* \to H^0(Z,\CS^*_{|Z})$.
Thus we obtain a rational map $H(X) \to \Gr(n-1,W^*) = \Gr(n+1,W)$.
More precisely, let
\begin{multline*}
\TTH(X) = \left\{(U,Z) \in \Gr(n+1,W)\times\Hilb_{d_n}(X)\
\left|
\begin{array}{l}
Z \in \Hilb_{d_n}(X\cap\Gr(2,U)),\\
\rank(\Lambda^2U^* \to H^0(Z,\Lambda^2\CS^*_{|Z})) \le c_n
\end{array}
\right\}\right..
\end{multline*}
We have projections $\TTH(X) \to \Gr(n+1,W)$ and $\TTH(X) \to H(X)$.

\begin{lemma}
The image of $\TTH(X)$ in $\Gr(n+1,W)$ coincides with $G(Y)$.
\end{lemma}
\begin{proof}
Let $(U,Z) \in \TTH(X)$. Since $Z \subset X\cap \Gr(2,U)$ we have commutative diagram
$$
\xymatrix{
V \ar[r] \ar@{-->} [dr] \ar@{..>}[drr]
& \Lambda^2W^* \ar[r] \ar[d] & \Lambda^2U^* \ar[d] \\
& H^0(X,\Lambda^2\CS^*_{|X}) \ar[r] & H^0(Z,\Lambda^2\CS^*_{|Z})
}
$$
Note that the dashed map vanishes by definition of $X$.
Therefore, the dotted arrow is also zero.
In other words, the composition
$$
V \to \Lambda^2U^* \to H^0(Z,\Lambda^2\CS^*_{|Z})
$$
is zero. On the other hand, the rank of the second map here is not greater than $c_n$ by definition of $\TTH(X)$.
Therefore the rank of the map $V \to \Lambda^2U^*$ is not greater than
$\dim\Lambda^2U^* - c_n = n(n+1)/2 - (n^2 - 3n + 4)/2 = 2n - 2$, so we see that $U \in G(Y)$ by definition
of $G(Y)$.

Vice versa, assume that $U \in G(Y)$. Then $X \cap \Gr(2,U) \subset \Gr(2,W)$
is the linear section of $\Gr(2,U)$ by the image of $V \subset \Lambda^2W^*$
in $\Lambda^2U^*$. By definition of $G(Y)$ the dimension of this image is
not greater than $2n - 2$, so $X \cap \Gr(2,U)$ is a linear section of the Grassmannian $\Gr(2,U) = \Gr(2,n+1)$
of codimension $\le 2n - 2$.
But since $\dim \Gr(2,n+1) = 2n - 2$
such an intersection contains not less than $\deg \Gr(2,n+1) = d_n$ points,
so we conclude that there exists a subscheme
$Z \subset X \cap \Gr(2,U)$ of length $d_n$.
Then $(U,Z) \in \TTH(X)$.
\end{proof}

Denote the maps from $\TTH(X)$ to $H(X)$ and $G(Y)$ by $\pi'$ and $\rho'$ respectively.
Consider the open subsets
$$
\begin{array}{l}
G'_0(Y) = \{ U \in \Gr(n+1,W)\ |\ \text{$X \cap \Gr(2,U)$ is zero-dimensional}\},\\
H_0(X) = \{ Z \in H(X)\ |\ \rank(W^* \to H^0(Z,\CS^*_{|Z})) = n+1 \}.
\end{array}
$$
We have a diagram
$$
\xymatrix{
& \TTH(X) \ar[dr]^{\pi'} \ar[dl]_{\rho'} \\
G'_0(Y) \subset G(Y) \qquad && \qquad H(X) \supset H_0(X)
}
$$
Also denote
$$
\begin{array}{ll}
\TTH'_0(X) = \rho'{}^{-1}(G'_0(Y)) = &\left\{\hbox to 0pt{$(U,Z)$\hss}\hphantom{Z \in \Hilb_{d_n}(X)}\ 
\left|\begin{array}{l}
Z = X\cap\Gr(2,U) \in \Hilb_{d_n},\\
\rank(\Lambda^2U^* \to H^0(Z,\Lambda^2\CS^*_{|Z})) \le c_n
\end{array}\right\}\right.,\\
\TTH_0(X) = \pi'{}^{-1}(H_0(X)) = &\left\{Z \in \Hilb_{d_n}(X)\ 
\left|\begin{array}{l}
\rank(W^* \to H^0(Z,\CS^*_{|Z})) = n+1,\\
\rank(\Lambda^2W^* \to H^0(Z,\Lambda^2\CS^*_{|Z})) \le c_n
\end{array}
\right\}\right..
\end{array}
$$

\begin{lemma}
If $n \ge 3$ then for a generic $X$ the set $\TTH'_0(X) \cap \TTH_0(X)$ is nonempty.
In particular, the sets $\TTH(X)$ and $H(X)$ have irreducible components
$\TTH_1(X)$ and $H_1(X)$ such that restriction of the maps $\rho'$ and $\pi'$
to $\TTH_1(X)$ are birational transformations $\TTH_1(X) \to G(Y)$ and $\TTH_1(X) \to H_1(X)$.
\end{lemma}
\begin{proof}
Once again, consider the universal versions of $\TTH(X)$, $\TTH_0(X)$ and $\TTH'_0(X)$:
\begin{multline*}
\tilde\CH_0,\tilde\CH'_0 \subset \tilde\CH = \{(V,U,Z) \in \Gr(2n,\Lambda^2W^*)\times\Gr(n+1,W)\times\Hilb_{d_n}(\Gr(2,W))\ |\\
|\
Z \subset \PP(V^\perp) \cap \Gr(2,U) \subset \Gr(2,W)\ \text{and}\
\rank(\Lambda^2U^* \to H^0(Z,\Lambda^2\CS^*_{|Z})) \le c_n
\}.
\end{multline*}
It is easy to see that $\tilde\CH_0 \cap \tilde\CH'_0 \subset \tilde\CH$ is a nonempty open subset
and its projection to $\Gr(2n,\Lambda^2W^*)$ is dominant. It follows that the generic fiber
of $\tilde\CH_0 \cap \tilde\CH'_0$ over $\Gr(2n,\Lambda^2W^*)$ is nonempty, so
$\TTH_0(X) \cap \TTH'_0(X)$ is nonempty for generic $X$. Define $\TTH_1(X)$ to be the closure
of $\TTH_0(X) \cap \TTH'_0(X)$ in $\TTH(X)$ and $H_1(X) = \pi'(\TTH_1(X))$.
Then the claim becomes obvious.
\end{proof}

We have proved:

\begin{theorem}\label{FY_bir_HX}
For generic $V \subset \Lambda^2W^*$, the three varieties $F(Y)$, $G(Y)$ and $H(X)$ are birational.
\end{theorem}

%
%

Since $H(X)$ is a moduli space on a Calabi--Yau manifold $X$ of dimension
$2n-4$, it has a natural $(2n-4)$-form.

\begin{conjecture}\label{4-forms_compatible}
The birational isomorphism of $F(Y)$ and $H(X)$ is compatible with $(2n-4)$-forms.
\end{conjecture}

We conclude this section by describing the above construction in cases $n=2$ and $n=3$ more explicitly.

If $n = 2$ then $\Pf(\CC^4) = \Gr(2,4)$, so $Y = \PP^3 \cap \Gr(2,4)$ is a quadric in $\PP^3$,
and $X = \Gr(2,4) \cap \PP^1$ is a pair of points. In this case $d_n = 1$, $c_n = 1$,
so $H(X) = X$ is a pair of points. However, the map $\pi'$ in this case is not birational,
it is a $\PP^1$-fibration in fact, so $\TTH(X)$ is a union of two $\PP^1$,
as well as $G(Y)$ and $F(Y)$. So, in case $n = 2$ our construction shows that
the variety of lines on a 2-dimensional quadric is $\PP^1 \sqcup \PP^1$.

If $n = 3$ then $Y = \PP^5 \cap \Pf(\CC^6)$ is a Pfaffian cubic fourfold,
and $X = \Gr(2,6) \cap \PP^8$ is a K3-surface. In this case $d = 2$, $c_n = 2$
and it is easy to see that all the maps $\pi,\rho,\rho',\pi'$ described above
are biregular isomorphisms. Moreover, the conditions
$\rank(W^* \to H^0(Z,\CS^*_{|Z})) \le 4,$
$\rank(\Lambda^2W^* \to H^0(Z,\Lambda^2\CS^*_{|Z})) \le 2$
defining $H(X) \subset \Hilb_2(X)$
are void since $Z$ is a length 2 subscheme in $X$,
so $H(X) = \Hilb_2(X)$ and our construction gives the classical
isomorphism between the Fano scheme $F(Y)$ of lines of $Y$ and $\Hilb_2(X)$.

\section{Exterior forms via the linkage class}\label{eflc}

We will start by defining the divisorial linkage class (see \cite{KM}, Section
3 for a more general notion of $k$-th linkage classes associated to any closed
embedding $i:Y\into W$). In what follows, $Y$ is a hypersurface (that is,
an effective Cartier divisor) in a nonsingular algebraic variety $W$.
We reserve the calligraphic letters $\FFF$, $\GGG$,\ldots\ for coherent
$\OOO_Y$-modules or objects of $\D^b(Y)$ and the block letters $F$, $G$, \ldots\
for coherent $\OOO_W$-modules or objects of $\D^b(W)$.

The restriction (or pullback) functor $i^*:\Coh (W)\rar \Coh (Y)$
has a left derived functor $Li^*:\D^b(W)\rar \D^b(Y)$. It can be described as follows:
let $\RRR(F)$  denote a locally free resolution of any $F\in \D^b(W)$. Then
$Li^*(F)$ is quasi-isomorphic to $i^*\RRR(F)$. The cohomology $h^kLi^*(F)$ of the complex
$Li^*(F)$ is denoted by $L_{-k}i^*(F)$ or $L^{k}i^*(F)$. Assume that
$F$ is a sheaf, that is a complex concentrated at grade 0.
Then the cohomology of $Li^*(F)$ can only appear at $k\leq 0$. More
exactly, we can write
$Li^*(F)=\OOO_Y\otimes_{f^{-1}\OOO_W}f^{-1}\RRR(F)$ and compute $L_ki^*(F)$
using the symmetry of the tensor product on its
arguments together with the fact that $\OOO_Y$, as an $\OOO_W$-module,
has a locally free resolution of the form $\RRR(i_*\OOO_Y)=[\OOO_W(-Y)\rar \OOO_W]$.
\textcolor{black}{We will apply this to the sheaves of the form $F=i_*\FFF$, where
$\FFF\in \Coh (Y)$.
We obtain that $L_0i^*(i_*\FFF)=\FFF$, $L_1i^*(i_*\FFF)=\FFF\otimes\OOO_W(-Y)|_Y
=\FFF\otimes\NNN^\dual_{Y/W}$,} where $\NNN^\dual_{Y/W}$ denotes
the conormal sheaf of $Y$ in $W$, and all the other cohomologies are zero.

If we consider $\RRR(i_*\OOO_Y)\otimes_{f^{-1}\OOO_W}f^{-1}F$
and $\OOO_Y\otimes_{f^{-1}\OOO_W}f^{-1}\RRR(F)$ as  complexes of
$\OOO_W$-modules, then they are quasi-isomorphic to each other and
to the direct sum of their cohomologies. But they are different
as complexes of $\OOO_Y$-modules: the first one is quasi-isomorphic
to the direct sum of its cohomologies, the second may be a nontrivial
extension. More exactly, let $\tau^\fatdot Li^*(F)$ be the canonical
filtration on $Li^*(F)$. Then $\tau^{-1} Li^*(F)$ and $\gr^0_\tau Li^*(F)$
are complexes concentrated at only one grade, so they are quasi-isomorphic
to their cohomologies:
$$
\tau^{-1} Li^*(F)\qis L_1i^*(F)[1],\ \ \ \
\gr^0_\tau Li^*(F)\qis L_0i^*(F).
$$
Thus the exact triple
\begin{equation}
  \label{triple}
0\rar \tau^{-1} Li^*(F)
\rar \tau^{0} Li^*(F)\rar \gr^0_\tau Li^*(F)\rar 0
\end{equation}
defines the extension class
\begin{equation}
  \label{eps_F}
  \eps_\FFF\in \Ext^1(L_0i^*(F),L_1i^*(F)[1])
=\Ext^2(L_0i^*(F),L_1i^*(F))=\textcolor{black}{
\Ext^2(\FFF,\FFF\otimes\NNN^\dual_{Y/W}).}
\end{equation}

\begin{definition} \textcolor{black}{Let $\FFF\in\Coh(Y)$
and $F=i_*\FFF$.}
The extension class \eqref{eps_F} of the exact triple \eqref{triple}
is called the divisorial
linkage class of $\FFF$ with respect to the embedding $i:Y\into W$.
To make explicit its dependence on the embedding, we will also denote it
by  $\eps^{Y/W}_\FFF$.
\end{definition}

\begin{remark}\label{4-term}
In practice, to compute $\eps_\FFF$, one can use a partial resolution
$0\rar G\rar E\rar F\rar 0$ of $F=i_*\FFF$ with $E$ locally free and $G$ torsion free
instead of the full resolution $\RRR(F)$. Tensoring by $\OOO_Y$,
we get the 4-term exact sequence
$$
0\rar i^*\TOR_1(i_*\OOO_Y,F) \rar i^*G
\rar i^*E\rar i^*F\rar 0,
$$
whose extension class in $\Ext^2(i^*F,i^*\TOR_1(i_*\OOO_Y,F))
=\textcolor{black}{\Ext^2(\FFF,\FFF\otimes\NNN^\dual_{Y/W}) }$ is $\eps_\FFF$.
\end{remark}

\begin{remark}\label{triangle}
In \cite{KM}, a more abstract definition of $\eps_\FFF$ is given,
valid for an object $\FFF$ of $\D^b(Y)$. It is a morphism
in the derived category included into a distinguished
triangle
$\xymatrix@1{Li^*i_*\CF \ar[r] & \CF \ar[r]^-{\epsilon_\CF} & \CF\otimes\NNN^\dual_{Y/W}[2]}$.
\end{remark}

The following example seems to be the easiest one for which
the linkage class is nonzero.

\begin{example}
Take for $i$ the Segre embedding $Y=\PP^1\times\PP^1\into W=\PP^3$
with image given by the equation $xw-yz=0$, and set $\FFF=\OOO_Y(1,0)$.
Then $F=i_*\FFF$ has a resolution $[\OOO_W(-1)^{\oplus 2}
\rar \OOO_W^{\oplus 2}]$ with the map given by the matrix
$\begin{pmatrix}x&y\\ z&w
\end{pmatrix}
$. It follows that $Li^*F\qis [\OOO_Y(-1,-1)^{\oplus 2}
\rar \OOO_Y^{\oplus 2}]$ with the map given by the same matrix,
$L_1i^*(F)=\OOO_Y(-1,-2)$, $L_0i^*(F)=\FFF=\OOO_Y(1,0)$,
and, by Remark \ref{4-term}, $\eps_{\FFF}\in H^2(Y,\OOO_Y(-2,-2))$ is the extension class
of the exact quadruple
$$
\arraycolsep=0.cm
\begin{array}{rcccccl}
0\rar\OOO_Y(-1,-2)&\xrightarrow{\phantom{xxxxx}}&\OOO_Y(-1,-1)^{\oplus 2}&
\xrightarrow{\phantom{xxxxxx}} &\OOO_Y^{\oplus 2}&\xrightarrow{\phantom{xxxxxxx}}&
\OOO_Y(1,0)\rar 0\\
&\begin{pmatrix}-y\\x
\end{pmatrix}&&\begin{pmatrix}x\ \ &y\\ z\ \ &w
\end{pmatrix}&&\begin{pmatrix}w\ \ &-y
\end{pmatrix}&
\end{array}
$$
\end{example}

We have $\eps_{\FFF}=0\IFF  Li^*F\qis\OOO_Y(1,0)\oplus
\OOO_Y(-1,-2)[1]\implies H^0(Y,\OOO_Y(-1,0)\otimes  Li^*F)=\CC$.
But from the resolution,
$$
H^0(Y,\OOO_Y(-1,0)\otimes  Li^*F)=
H^0(Y,[\OOO_Y(-2,-1)^{\oplus 2}\rar \OOO_Y(-1,0)^{\oplus 2}])=0,
$$
hence $\eps_{\FFF}\neq 0$.

Now we go over to the situation when $Y=Y_n$ is a hypersurface
of degree $n$ in $W=\PP^{2n-1}$ and $\FFF$ is a sheaf representing
a point $m=[\FFF]$ of some connected component $M$ of the smooth locus
of the moduli space of stable (or simple) sheaves
on $Y$. There is a canonical isomorphism $T_mM=
\Ext^1(\FFF,\FFF)$, and formula \eqref{p-form} defines
a $p$-linear form $\alp_p(m)$ on $T_mM$ with $p=2n-4$. It depends on the
choice of the isomorphism $\NNN^\dual_{Y/W}\simeq
\omega_Y$, where $\omega_Y$ is the dualizing sheaf of $Y$.
Fixing once and forever such an isomorphism, the forms  $\alp_p(m)$
fit to a well-defined cross-section $\alp_p$ of the vector bundle
$\wedge^p\TTT^\dual M=\Omega^{2n-4}_{M}$.
One way to see that it is a regular section is to relativize
the definition of $\eps_\FFF$
and formula \eqref{p-form} in flat families of sheaves. We will use
another approach, consisting in relating $\eps_\FFF$ to the Atiyah
class of $\FFF$.

We will use the Atiyah class of torsion sheaves on $Y$ with support
$Z$ which is a locally complete intersection subscheme in $Y$.
Let $i:Z\into Y$ be
the natural embedding. For any $\FFF,\GGG
\in \D^b(Z)$, we have a
spectral sequence
$$
E_2^{pq}=\Ext^p(L_qi^*i_*\CF,\GGG)\implies \Ext^{p+q}(i_*\CF,i_*\GGG),
$$
and a canonical isomorphism $L_qi^*i_*\CF\isom{\rm can} \CF\otimes \wedge^q\CN^\dual_{Z/Y}$.
This provides a natural map $\Ext^{1}(i_*\CF,i_*\GGG)
\rar \Hom (\CF\otimes\CN_{Y/W}^\dual,\GGG)$.

Assume now that, moreover,
$\FFF$ is a locally free sheaf on $Z$. Then there is a
canonical isomorphism $
i^*\EXT^q(i_*\CF,i_*\GGG)=\HOM(\CF\otimes \wedge^q\CN^\dual_{Z/Y},\GGG)$,
and we obtain another natural map between the same objects:
\begin{multline*}\Ext^{1}(i_*\CF,i_*\GGG)\rar
H^0(Y,\EXT^{1}(i_*\CF,i_*\GGG))=H^0(Z,i^*\EXT^{1}(i_*\CF,i_*\GGG))=\\
=H^0(Z, \HOM (\CF\otimes\CN_{Z/Y}^\dual,\GGG))=\Hom (\CF\otimes\CN_{Z/Y}^\dual,\GGG),
\end{multline*}
where the first arrow comes from the local-to-global spectral sequence for Exts.
One can prove that these two maps coincide.

\begin{theorem} \label{kappat}
Let
$Z$ be a locally complete intersection subscheme in $Y$ and
$i:Z\into Y$ the natural embedding. Let $\CF \in \D^b(Z)$.
Then the image of the Atiyah class
$\At_{i_*\CF} \in \Ext^1(i_*\CF,i_*\CF\otimes\Omega^1_Y)$ in
$\Hom(\CF\otimes\CN_{Z/Y}^\dual,\CF\otimes\Omega^1_Y|_Z)$
coincides with $1_\CF\otimes\kappa_{Z/Y}$,
where $\kappa_{Z/Y}$ denotes the natural
map of sheaves $\CN^\dual_{Z/Y}\rar \Omega^1_Y|_Z$.
\end{theorem}

\begin{proof}
See \cite{KM}, Theorem 3.2 (iii).
\end{proof}

When $Y$ is a hypersurface in a smooth variety $W$, consider the conormal bundle  sequence
\begin{equation}\label{cnbs}
\xymatrix@1{0 \ar[r]& \CN^\dual_{Y/W} \ar[r]^{\ \kappa_{Y/W}\ }&
  \Omega^1_W|_Y \ar[r]^{\ \rho_{Y/W}\ }& \Omega^1_Y \ar[r]& 0}.
\end{equation}
Denote by $\nu_{Y/W}$ the extension class of~\eqref{cnbs} in  $\Ext^1(\Omega^1_Y,\CN^\dual_{Y/W})$.

\begin{theorem}\label{epsat}
Let $Y$ be a hypersurface in a smooth variety $W$, $i:Y \to W$ the natural embedding, and  $\CF \in \D^b(Y)$. Then
the linkage class $\epsilon^{Y/W}_\CF \in \Ext^2(\CF,\CF\otimes\CN^\dual_{Y/W})$
factors through the Atiyah class $\At_\CF \in \Ext^1(\CF,\CF\otimes\Omega^1_Y)$:
$$\epsilon^{Y/W}_\CF = (1_\CF\otimes\nu_{Y/W})\circ\At_\CF.$$
\end{theorem}

\begin{proof}
See \cite{KM}, Theorem 3.2 (i) \textcolor{black}{or \cite{HT}, Proposition~3.1.}
\end{proof}

In the next theorem, we deal with exterior forms on a moduli space $M$
of sheaves on a projective variety $Y$.
Speaking about exterior forms on $M$, we will always mean
that $M$ is smooth.
On the other hand, it is not important to us whether
$M$ is separated or not. Thus we will take for $M$ an open subset of
the smooth locus of either one of the following two moduli spaces:
either the moduli space of $H$-stable sheaves
on $Y$ in the sense of the definition of Simpson \cite{S} for some
ample divisor class $H$ on $Y$, or
the moduli space of simple sheaves on $Y$ as defined in \cite{AK2}.
The former is a quasi-projective scheme over $k$,
and the latter is a possibly non-separated algebraic space.

\begin{corollary}\label{closed}
Let $Y=Y_n$ be a smooth hypersurface
of degree $n$ in $W=\PP^{2n-1}$ and $M$ a connected component of the smooth locus
of the moduli space of stable or simple sheaves
on $Y$. Let us fix an isomorphism $\NNN^\dual_{Y/W}\simeq
\omega_Y$. Then there exists a closed regular $p$-form $\alp_p
\in H^0(M, \Omega^p_M)$, where $p=2n-4$, such that its value $\alp_p(m)$
at any point $m=[\FFF]\in M$ represented by a sheaf $\FFF$ is the composition
\begin{multline}\label{p-form}
\underbrace{\Ext^1(\FFF,\FFF)\times\ldots\times
\Ext^1(\FFF,\FFF)}_{\mbox{$2n-4$ times}}
\xrightarrow{\mbox{\scriptsize Yoneda}}\Ext^{2n-4}(\FFF,\FFF)
\xrightarrow{\eps_\FFF}\Ext^{2n-2}(\FFF,\FFF\otimes
\NNN^\dual_{Y_n/\PP^{2n-1}})\\ \simeq
\Ext^{2n-2}(\FFF,\FFF\otimes
\omega_{Y_n})
\xrightarrow{\Tr}H^{2n-2}(Y_n, \omega_{Y_n})=\CC.
\end{multline}
If $n=4$ and $M=F(Y)$, $\alp_p$ is proportional to the $4$-form on $F(Y)$
defined in Corollary \ref{4_form}.
\end{corollary}

\begin{proof}
As in the proof of Theorem 2.2 of \cite{KM}, we can shrink $M$ to the
biholomorphic image of a polydisk in it, choose a universal
sheaf $\FFFF$ over $M\times Y$ and represent
$(p+1)\alp_p$ as the K\"unneth component $\gamma^{2n-4,2n-2}\in
H^0(M,\Omega_M^{2n-4})\otimes H^{2n-2}(Y,\Omega^{2n-2}_Y)$
of
$$\gamma=\Tr(\At_\FFFF^{\wedge (2n-3)})\wedge \nu_{Y/W}\in
H^{2n-2}(M\times Y, \Omega^{4n-6}_{M\times Y}),
$$
where $\nu_{Y/W}$ is viewed here as an element of $H^{1}(Y,\Omega^{2n-3}_Y)\simeq
\Ext^1(\Omega^1_Y,\CN^\dual_{Y/W})$. Then $\gamma$ is de Rham closed by
\cite{HL}, Sect. 10.1.6. In the same way as in loc. cit., this fact together with
the projectivity of $Y$ implies that all the K\"unneth components
of $\gamma$ are $d_M$-closed,
where $d_M$ denotes the de Rham differential on $M$.
\end{proof}

To conclude this section, we will prove a nonvanishing theorem for
$\eps_\FFF$ which is a direct generalization of Proposition 4.1 of \cite{KM}
from $n=3$ to all $n\geq 3$. We use here a more general definition
of the linkage class, referred to in Remark \ref{triangle}.

\begin{proposition}\label{lc-on-exts}
Assume that $\CF,\CG \in \CCC_n$, that is $H^p(Y,\CF(-k))=H^p(Y,\CG(-k))=0$
for all $p\in \ZZ$ and $k=0,1,\ldots ,n-1$.
Then the multiplication by the linkage class
$\epsilon_\CG \in \Ext^2(\CG,\CG(-n))$ induces an isomorphism
$\Ext^p(\CF,\CG) \cong \Ext^{p+2}(\CF,\CG(-n))$ for all $p\in\ZZ$.
\end{proposition}
\begin{proof}
Consider the Beilinson spectral sequence for $i_*\CG$ (see \cite{Bei,OSS})
$$
E_1^{-p,q} = H^q(\PP^{2n-1},i_*\CG(-p))\otimes \Omega_{\PP^{2n-1}}^p(p)
\Longrightarrow i_*\CG ,
$$
where $i:Y\hookrightarrow \PP^{2n-1}$ is the natural embedding.
By the assumption on $\GGG$,
we have $E_1^{0,q} = E_1^{-1,q} =\ldots = E_1^{-n+1,q} = 0$ for all $q$. Hence
the derived pullback $Li^*i_*\CG$ is contained in the triangulated
subcategory of $\D^b(\Coh(Y))$ generated by $i^*\Omega_{\PP^{2n-1}}^n(n)$,
$i^*\Omega_{\PP^{2n-1}}^{n+1}(n+1)$, \ldots\ ,
$i^*\Omega_{\PP^{2n-1}}^{2n-1}({2n-1})$. Let $V=\CC^{2n}$, so that $\PP^{2n-1}=\PP(V)$.
The Euler exact sequence
$$
0\lra \OOO_{\PP^{2n-1}}\lra V\otimes\OOO_{\PP^{2n-1}}(1)
\lra \TTT_{\PP^{2n-1}}\lra 0
$$
and the isomorphism $\Omega_{\PP^{2n-1}}^{2n-1-k}
\simeq \wedge^k\TTT_{\PP^{2n-1}}(-2n)$
lead to the following resolutions for the sheaves $\Omega^{2n-1-k}(2n-1-k)$
($k=0,1,\ldots, n-1$):
$$
0\rar \OOO(-k-1)\rar\ \ \ldots\ \ \rar
\wedge^{k-1}V\otimes \OOO(-2)\rar
\wedge^kV\otimes \OOO(-1)\rar
\Omega^{2n-1-k}(2n-1-k)\rar 0\ .
$$
Hence
this subcategory coincides with the subcategory of $\D^b(\Coh(Y))$ generated
by $\CO_Y(-1)$, \ldots\ , $\CO_Y(-n)$. By Serre duality on $Y$ and by the hypothesis on
$\FFF$, we obtain:
$$
\Ext^p(\CF,\CO_Y(-k-1)) \cong H^{2n-1-p}(\CF(-k))^\dual = 0
$$
for all $k=0,1,\ldots, n-1$.
Hence $\Ext^\bullet(\CF,Li^*i_*\CG) = 0$.
It remains to note that we have a distinguished triangle
$\xymatrix@1{Li^*i_*\CG \ar[r] & \CG \ar[r]^-{\epsilon_\CG} & \CG(-n)[2]}$.
Applying the functor $\Hom(\CF,-)$, we deduce the proposition.
\end{proof}

For $n=3$, this nonvanishing property implies
the nondegeneracy of $\alp_{2n-4}$ on moduli spaces $M$ that parameterize
sheaves from $\CCC_n$. The same argument implies only partial results for bigger $n$.
For example:

\begin{corollary}
Let $n=4$, so that $Y$ is a quartic in $\PP^7$, and assume that a moduli
space $M$ parameterizes sheaves from $\CCC_n$. Then the $2$-rank of $\alp_4$
at a point $m\in M$ representing a sheaf $\FFF$
coincides with the dimension of the image of the Yoneda coupling
$$
\Ext^1(\FFF,\FFF)\times \Ext^1(\FFF,\FFF)\lra \Ext^2(\FFF,\FFF).
$$
\end{corollary}

\begin{proof} Similar to the proof of Theorem 4.3 in \cite{KM}.
\end{proof}

\section{An explicit calculation}\label{ec}

Throughout this section, $Y$ will denote a generic
quartic hypersurface in $\PP^7$ (a quartic 6-fold).
Let $F=F(Y)$ be the Fano scheme of $Y$, that is the Hilbert scheme of lines in $Y$.
We have seen that $F$ is a smooth connected projective variety of dimension $7$.
We will denote by $\{\ell\}$ the point of $F(Y)$ representing
a line $\ell\subset Y$.

The constructions of Corollaries \ref{4_form} and \ref{closed} provide a closed 4-form
$\alp=\alp_4\in H^0(F(Y),\Omega^4_{F(Y)})$. We will produce
an explicit formula for its value $\alp (m)$ at a point $m=\{\ell\}$.
\textcolor{black}{
This calculation will allow us to prove in Section \ref{4-forms} that
$\alp_4$ is a minimally degenerate 4-form on} \textcolor{black}{$F(Y_4)$.}
\textcolor{black}{Moreover, it
indicates a general pattern of such calculation for $Y_n$ ($n>4$), which
is only notationally harder than that for $Y_4$. }

\begin{proposition}
Let $Y$ be a generic quartic hypersurface in $\PP^7$. Then
there is an algebraic subset $R(Y)$ in $F(Y)$ of codimension $\geq 3$
such that
$\NNN_{\ell/Y}\simeq \OOO(1)^{\oplus 2}\oplus\OOO^{\oplus 3}$ for all
$\{\ell\}\in F(Y)\setminus R(Y)$, and
$\NNN_{\ell/Y}\simeq \OOO(1)^{\oplus 3}\oplus\OOO\oplus\OOO(-1)$
if $\{\ell\}\in R(Y)$.
\end{proposition}

\begin{proof}
Denote by $V$ the vector space $\CC^8$ whose projectivization contains $Y$.
So $\PP^7=\PP(V)$, and $Y=Y_f$ is defined by an octal quartic form $f\in S^4V^*$.

The possible types of the normal bundle follow easily from the two exact triples
$$
0\rar \TTT_\ell\rar\TTT_Y|_\ell\rar\NNN_{\ell/Y}\rar 0\ ,\ \ \
0\rar \TTT_Y|_\ell \rar\TTT_{\PP^7}|_\ell\rar\NNN_{Y/\PP^7}|_\ell\rar 0.
$$
Indeed, the tangent space $T_{\{\ell\}}F(Y)$ is of dimension
$h^0(\NNN_{\ell/Y})$, so the smoothness of $F(Y)$ implies that
$h^1(\NNN_{\ell/Y})=0$ for all $\ell$. Hence the Grothendieck
splitting of $\NNN_{\ell/Y}$ has no summands $\OOO(a)$ with
$a<-1$. As $\TTT_{\PP^7}|_\ell\simeq \OOO(2)\oplus\OOO(1)^{\oplus 6}$
and $\TTT_Y|_\ell$ has an injective map to it, $\TTT_Y|_\ell$
is the sum of sheaves $\OOO(a_i)$ with $2\geq a_1\geq 1\geq a_2\geq\cdots
\geq a_6\geq -1$, and $\sum a_i=\deg \TTT_Y|_\ell=4$. Moreover, $a_1=2$
because $\TTT_\ell\simeq \OOO(2)$ has a nonzero map to $\TTT_Y|_\ell$. This leaves
only two possible choices $(a_i)=
(2,1^2,0^3)$ or $(2,1^3,0,-1)$. Splitting off $\OOO(a_1)\simeq \TTT_\ell$,
we obtain the two possible normal bundles $\NNN_{\ell/Y}$.

Thus
the lines in $Y$ are of two types. The summands $\OOO(1)$ of $\NNN_{\ell/Y}$
span the subbundle corresponding to infinitesimal deformations of $\ell$
inside the projective subspace $\bigcap_{P\in\ell}T_PY$, the intersection
of the tangent hyperplanes to $Y$ at all the points of $\ell$.
This is a cubic pencil of hyperplanes, and generically there are 4
linearly independent ones. Hence $\bigcap_{P\in\ell}T_PY\simeq \PP^3$
for $\ell$ in an open subset of $F(Y)$. These are lines of the first type.
Their complement $R(Y)$ in $F(Y)$ is the set of lines of the second type;
for them $\bigcap_{P\in\ell}T_PY\simeq \PP^4$.
It remains to estimate the codimension of $R(Y)$.

Let $g=g_f:Y\rar Y^*\subset \PP^{7*}$, $P\mapsto T_PY$ be the Gauss map,
given in homogeneous coordinates by $(x_0:\cdots:x_7)
\mapsto (\partial f/\partial x_0:\cdots:\partial f/\partial x_7)$,
where $f$ is the quartic form defining $Y$. As the partial
derivatives $\partial f/\partial x_i$ have no common zero on $Y$,
the Gauss map is finite.
Its restriction to a line $\ell$ in $Y$ is given by a cubic pencil
without fixed points, so the image $g(\ell)$ is either
a cubic rational curve, or a line. In the latter case, $\bigcap_{P\in\ell}T_PY$
would be 5-dimensional, which is impossible by the above. Hence $g(\ell)$
is always a rational cubic, and $\{\ell\}\in R(Y)$ if and only if $g(\ell)$ has
exactly one singular point. Thus if we define the algebraic set
$$
\tilde R(Y)=\{(\ell, P)\in F(Y)\times Y\ \mid \ P\in\ell,\ g(P)\in\Sing g(\ell)\},
$$
then the projection $\pr_1:\tilde R(Y)\rar F(Y)$ is finite of degree 2
and has $R(Y)$ as its image, so $\dim R(Y)=\dim \tilde R(Y)$. We will deduce the dimension of $\tilde R(Y)$ by
a standard dimension count.

Consider two incidence varieties:
$$
I_1=\{ (P,Q,f)\in\PP(V)\times\PP(V)\times\PP(S^4V^*)\mid
P\neq Q,\ \nabla_Pf\sim\nabla_Qf,\ \overline{PQ}\subset Y_f\},
$$
where
$\nabla_Pf=(\partial f/\partial x_0,\cdots,\partial f/\partial x_7)|_P$
is the gradient of $f$ at $P$, the sign $\sim$ stands for proportionality, and
$\overline{PQ}$ denotes the line passing through $P,Q$, and
$$
I_2=\{ ([P,v],f)\in \PP(\TTT_{\PP^7})\times\PP(S^4V^*)\mid
P\in Y_f,\ v\in T_PY_f,\ v\neq 0,\ (\nabla_Pf,v)=0,\ \overline{Pv}\subset Y_f\},
$$
where $[P,v]$ is the proportionality class of $v$ considered as a
point of $\PP(T_P\PP^7)$, and $\overline{Pv}$ is the line through
$P$ in the direction of $v$.

The proportionality $\nabla_Pf\sim\nabla_Qf$ can be interpreted
as the coincidence of the tangent spaces $T_PY_f=T_QY_f$ as soon as both
gradients are nonzero.
The part of $I_1$ for which $Y_f$ is nonsingular
parametrizes all the triples $(P,Q,f)$ for which $g_f(\overline{PQ})$ has a node
at $g_f(P)=g_f(Q)$. Similarly,
the part of  $I_2$ for which $Y_f$ is nonsingular
parametrizes the pairs $([P,v],f)$ for which
$g_f(\overline{Pv})$ has a cusp at $g_f(P)$. Thus, assuming $Y=Y_f$
nonsingular, we can represent $\tilde R(Y)$ as the union of two
algebraic sets $\pi_{i}^{-1}(f)$, where $\pi_{i}:I_i\rar \PP(S^4V^*)$
($i=1,2$) is the natural projection.

Looking at the other projection $\pr_{12}:I_1\rar \PP(V)\times\PP(V)
\setminus\mbox{(diagonal)}$, we find that all of its fibers are
isomorphic to each other and irreducible
of dimension 319, hence $I_1$ is irreducible of dimension 333. As $\dim
\PP(S^4V^*)=329$, we conclude that $\pi_{1}^{-1}(f)$ is either empty,
or is of dimension $333-329=4$ for generic $f$. Similarly, $I_2$
is irreducible and $\dim I_2=332$, so $\pi_{2}^{-1}(f)$ is either empty
or is of dimension $3$ for generic $f$. Hence $\dim \tilde R(Y)\leq 4$
for generic $f$, as was to be proved.
\end{proof}

Fix now a generic quartic $Y\subset \PP^7=\PP (V)$ and a line
$\ell$ of first type in $Y$. Choose homogeneous coordinates in $\PP^7$
in such a way that $\ell=\{ x_0=\cdots =x_5=0\}$. Then the equation $f$
of $Y$ can be written in the form
$$f=x_0f_0(x_0,\ldots, x_7)+\cdots +x_5f_5(x_0,\ldots, x_7),$$
where the $f_i$ are cubic forms.
Denote by $\bar{f}_i=\bar{f}_i(x_6,x_7)$ the restriction of
$f_i$ to $\ell$. The fact that $\ell$ is of the first type
implies that the 6 forms $\bar{f}_i=\frac{\partial f}{\partial x_i}|_\ell$
generate the whole 4-dimensional
vector space of binary cubics, and by a linear change of coordinates
$x_0,\ldots, x_5$, we can arrange the things so that
$$
\bar f_0=x_6^3,\
\bar f_1=x_6^2x_7,\
\bar f_2=x_6x^2_7,\
\bar f_3=x_7^3, \
\bar f_4=\bar f_5=0.
$$
Now we will construct an explicit Grothendieck splitting of $\TTT_Y|_\ell$
and $\NNN_{\ell/Y}$. We have $\PP (V)=V^0/\CC^*$,
where $V^0=V\setminus\{0\}$. The symbols $\partial/\partial x_i$,
$dx_i$ have a natural meaning as vector fields, respectively 1-forms
on $V$ or $V^0$. A rational section of $\TTT_{\PP^7}(k)$ can be
represented in the form $\sum \phi_i\partial/\partial x_i\mod E$,
where $\phi_i$ are rational homogeneous functions in $x_j$ of
degree $k+1$, and $E=\sum x_i\partial/\partial x_i$ is the Euler
vector field. Let us denote $\partial/\partial x_i\mod E$ by
$\partial_i$. Then $\partial_0,\ldots,\partial_7$ form a basis
of $H^0(\TTT_{\PP^7}(-1))$, and $H^0(\TTT_{\PP^7})$ is generated
by $x_i\partial_j$, $0\leq i,j\leq 7$, with a single linear
relation $\sum x_i\partial_i=0$.

A rational vector field $\sum \phi_i\partial_i$ is tangent to $Y$
if and only if $\sum\phi_i\partial f/\partial x_i|_Y=0$. Similarly,
an expression $\sum\psi_idx_i$, where $\psi_i$ are rational
homogeneous functions in $x_j$ of degree $k-1$, represents
a rational section of $\Omega^1_{\PP^7}(k)$ whenever
$\sum x_i\psi_i=0$, and then
$\sum\psi_idx_i|_Y\mod df$ represents
a rational section of $\Omega^1_{Y}(k)$. We will use the overbar
to denote the restriction of sections of $\Omega^1_{\PP^7}(k)$ or
$\TTT_{\PP^7}(k)$ to sections of $\Omega^1_{Y}(k)|_\ell$ or $\TTT_{Y}(k)|_\ell$.
With this notation, we have the following Grothendieck splitting:
\begin{equation}\label{spl}
\TTT_Y |_\ell=\OOO_\ell(2)e_0\oplus \OOO_\ell(1) e_1\oplus\OOO_\ell(1) e_2
\oplus\OOO_\ell e_3\oplus\OOO_\ell e_4\oplus\OOO_\ell e_5,
\end{equation}
where the $e_i$ are the following sections of appropriate twists of $\TTT_Y |_\ell$:
\begin{multline}\label{basis}
e_0=x_6^{-1}\bar\partial_7=-x_7^{-1}\bar\partial_6,\
e_1=\bar\partial_4,\ e_2=\bar\partial_5,\\
e_3=x_6\bar\partial_1-x_7\bar\partial_{0},\ e_4=x_6\bar\partial_2-x_7\bar\partial_{1},\  e_5=x_6\bar\partial_3-x_7\bar\partial_{2}.
\end{multline}

Thus $e_0\in H^0(\TTT_{Y}(-2)|_\ell)$, $e_1$ and $e_2$ are elements
of $H^0(\TTT_{Y}(-1)|_\ell)$, and the remaining $e_i$ are sections
of $\TTT_{Y}|_\ell$. We can identify $\NNN_{\ell/Y}$ as the sum
of the last five summands. Denote by $\check e_i$ the dual basis of
$H^0(\ell, (\Omega^1_Y|_\ell)\otimes k(\ell))$, so that
$(\check e_i,e_j)=\delta_{ij}$. Then
\begin{multline}\label{db}
\check e_0=\bar{x_6dx_7-x_7dx_6}\in H^0(\Omega^1_Y(2)|_\ell),\\
\check e_1=d\bar x_4,\ \check e_2=d\bar x_5\in H^0(\Omega^1_Y(1)|_\ell),\ \
H^0(\Omega^1_Y|_\ell)=\langle \check e_3,\check e_4, \check e_5\rangle,\
\check e_3=-\frac{d\bar x_0}{x_7},\\
\check e_4=\frac{1}{2}\left(-\frac{x_6}{x_7^2}d\bar x_0-\frac{1}{x_7}d\bar x_1
+\frac{1}{x_6}d\bar x_2+\frac{x_7}{x_6^2}d\bar x_3\right),\ \ \check e_5=
\frac{d\bar x_3}{x_6}
.\end{multline}
Here $d\bar x_i$ are well-defined for $i=0,\ldots,5$ by the
formula
\begin{equation}\label{dx-bar}
d\bar x_i=\bar{dx_i-\frac{x_i}{x_j}dx_j}=\left.\left(
dx_i-\frac{x_i}{x_j}dx_j\right)\right|_{\textstyle \ell}\mod
df|_{\textstyle \ell}\ \ (j=6\ \mbox{or}\ 7).
\end{equation}
The Grothendieck splitting of $\Omega^1_Y|_\ell$ has
the form
\begin{equation}\label{dspl}
\Omega^1_Y|_\ell=\OOO_\ell(-2)\check e_0\oplus \OOO_\ell(-1) \check e_1\oplus\OOO_\ell(-1) \check e_2
\oplus\OOO_\ell \check e_3\oplus\OOO_\ell \check e_4\oplus\OOO_\ell \check e_5
\end{equation}
and the last four summands give a Grothendieck splitting of
$\NNN_{\ell/Y}^\dual$.

Let us now compute the composition (\ref{p-form}) with $p=4$
on four sections $\xi_i\in \Ext^1(\OOO_\ell,\OOO_\ell)=H^0(\NNN_{\ell/Y})$,
$i=1,\ldots,4$,
\begin{equation}\label{xii}
\xi_i=(a_{i0}x_6+a_{i1}x_7)e_1+(b_{i0}x_6+b_{i1}x_7)e_2+c_{i1}e_3+c_{i2}e_4+c_{i3}e_5,\end{equation}
where $a_{ij},b_{ij},c_{ij}$ are elements of $k$.

\begin{theorem}\label{determinants}
In the above notation, the $4$-form $\alp_4(m)$ defined by \eqref{p-form}
at a point $m=[\OOO_\ell]$ of the Fano scheme $F(Y)$
is given up to a constant factor by the formula
\begin{multline}\label{4dets}
\alp_4(m)(\xi_1,\ldots,\xi_4)=\begin{vmatrix}
     a_{10}&b_{10}&c_{12}&c_{13}\\
     \vdots&\vdots&\vdots&\vdots\\
     a_{40}&b_{40}&c_{42}&c_{43}
     \end{vmatrix} \ \ - \ \
\begin{vmatrix}
a_{10}&b_{11}&c_{11}&c_{13}\\
\vdots&\vdots&\vdots&\vdots\\
a_{40}&b_{41}&c_{41}&c_{43}
\end{vmatrix} \\ - \ \
      \begin{vmatrix}
      a_{11}&b_{10}&c_{11}&c_{13}\\
      \vdots&\vdots&\vdots&\vdots\\
      a_{41}&b_{40}&c_{41}&c_{43}
      \end{vmatrix} \ \ + \ \
\begin{vmatrix}
a_{11}&b_{11}&c_{11}&c_{12}\\
\vdots&\vdots&\vdots&\vdots\\
a_{41}&b_{41}&c_{41}&c_{42}
\end{vmatrix}\ .
\end{multline}
\end{theorem}

\begin{proof}
By Theorem \ref{epsat}, $\alp_4(m)$
is the composition of the maps in the upper line
of the diagram
$$
\xymatrix@C=22pt{ \Ext^1(\OOO_\ell,\OOO_\ell)^4
\ar[r]^-{\rm Yoneda} \ar[d]_(.47){\simvert}&
\Ext^4(\OOO_\ell,\OOO_\ell) \ar[r]^-{\At_{\OOO_\ell}} \ar[d]_(.47){\simvert} &
 \Ext^5(\OOO_\ell,\OOO_\ell\otimes\Omega^1_Y) \ar[d]_(.47){\simvert}
                                                                 \ar[r]^-{\nu} &
  \Ext^6(\OOO_\ell,\OOO_\ell\otimes\omega_Y)\ar[d]_(.47){\simvert}\ar[r]^-{\Tr} &
                                                                    k \ar@{=}[d] \\
H^0(\NNN_{\ell/Y})^4  \ar[r]^-{\wedge} &
H^0(\wedge^4\NNN_{\ell/Y})  \ar[r]^-{\kappa} &
   H^0(\wedge^5\NNN_{\ell/Y}\otimes\Omega^1_Y) \ar[r]^-{\nu} &
                                H^1(\omega_\ell)\ar[r]^-{\Tr} &
                                                                k
}\ \
$$
Here we are using the canonical isomorphisms $\EXT^q(\OOO_\ell,\OOO_\ell)
=\wedge^q\NNN_{\ell/Y}$, and the vanishing $h^i(\wedge^q\NNN_{\ell/Y})
=h^i(\wedge^4\NNN_{\ell/Y}\otimes\Omega^1_Y)=0$ for $i>0$
implies that the vertical maps, coming from the local-to-global
spectral sequence,
are also isomorphisms. Further, by \cite{KM}, Lemma 1.3.2,
the Yoneda coupling on the Ext-sheaves
corresponds to the wedge product under these isomorphisms. The symbols
$\kappa$, $\nu$ were defined in Theorems \ref{kappat}, \ref{epsat}. After
applying $\nu\in\Ext^1(\Omega^1_Y,\NNN_{Y/\PP^7}^\dual)$ in both lines, we also use
the isomorphisms $\NNN_{Y/\PP^7}^\dual\simeq\omega_Y$,
$\wedge^5\NNN_{\ell/Y}\otimes\NNN_{Y/\PP^7}^\dual\simeq\omega_\ell\:$;
the constant factor mentioned in the statement of the theorem
depends only on the choice of these isomorphisms.

Thus we can calculate $\alp_4(m)$ in following the bottom line of the
digram, in which the Yoneda product turns into the wedge one.
We have $\kappa=e_1\otimes \check e_1+\cdots +e_5\otimes \check e_5$,
and the image of $(\xi_1,\ldots,\xi_4)$ in
$H^0(\wedge^5\NNN_{\ell/Y}\otimes\Omega^1_Y)$ is
$$
\psi=\xi_1\wedge \xi_2\wedge \xi_3\wedge \xi_4\wedge \kappa=
e_1\wedge\ldots \wedge e_5\otimes
\begin{vmatrix}
a_{10}x_6+a_{11}x_7&b_{10}x_6+b_{11}x_7&c_{11}&c_{12}&c_{13}\\
     \vdots&\vdots&\vdots&\vdots\\
 a_{40}x_6+a_{41}x_7&b_{40}x_6+b_{41}x_7&c_{41}&c_{42}&c_{43}\\
 \check e_1&\check e_2&\check e_3&\check e_4&\check e_5&
\end{vmatrix}\ ,
$$
where the determinant is understood as its row expansion with respect to
the last row. We can equally omit the factor
$e_1\wedge \ldots \wedge e_5$,
for it is a nowhere vanishing section
of $\det\NNN_{\ell/Y}(-2)\simeq \OOO_\ell$. Then $\psi$ is just the
section of $\Omega^1|_\ell(2)$ given by the determinant in the last
formula.

Next we have to couple $\psi$ with $\nu$. The Ext-group to which
$\nu$ belongs can be represented as $H^1(\TTT_Y(-4))$ or else $H^1(\Omega^5_Y)$.
As we know from  Proposition \ref{Hodge_num},
$h^{51}(Y)=1$.
Thus we have to determine the image of the 1-dimensional
$H^1(\Omega^5_Y)$ inside the 14-dimensional $H^1(\Omega^5_Y|_\ell)$.
Keeping in mind that $\Omega^5_Y\simeq \TTT_Y(-4)$, we use the natural exact
triple $$0\rar \TTT_Y(-4)\rar \TTT_{\PP^7}(-4)|_Y\rar \OOO_Y\rar 0.$$
Then we see that
$H^1(\TTT_Y(-4))=\beta_Y (H^0 (\OOO_Y))$, where
$\beta_Y$ is the Bockstein homomorphism, and the commutativity
of the maps of the exact triple with the restriction to $\ell$
implies that
$$\im \left(H^1(\TTT_Y(-4))\xrightarrow{\rm restriction}
H^1(\TTT_Y(-4)|_\ell)\right)\ =\ \im \left(H^0 (\OOO_\ell)\xrightarrow{
\beta_\ell}
H^1(\TTT_Y(-4)|_\ell)\right),$$
where $\beta_\ell$ is the Bockstein homomorphism of the restricted
exact triple
$$
0\lra \TTT_Y(-4)|_\ell\lra \TTT_{\PP^7}(-4)|_\ell\lra \OOO_\ell\lra 0.
$$
The surjection in this exact triple is induced by
$\bar \partial_i\mapsto (\partial f/\partial x_i)|_\ell$,
so that $x_6^{-3}\bar \partial_0\mapsto 1$ and
$x_7^{-3}\bar \partial_3\mapsto 1$. Taking the
standard covering of $\ell=\PP^1$ by the open affine sets
$U_6 =\{x_6\neq 0\}$, $U_7 =\{x_7\neq 0\}$,
we get a \v{C}ech representative for $\beta_\ell (1)\in H^1(\TTT_Y(-4)|_\ell)$
of the form
\begin{equation}\label{bock}
x_7^{-3}\bar \partial_3-x_6^{-3}\bar \partial_0=
\frac{x_7^2e_3+x_6x_7e_4+x_6^2e_5}{x_6^3x_7^3}
\in \Gamma(U_6\cap U_7, \TTT_Y(-4)|_\ell).
\end{equation}
The wanted quantity $\alp_4(m)(\xi_1,\ldots,\xi_4)$ is by construction $\Tr(\nu\wedge\psi)$ when $\nu$ is considered as an element
of $H^1(\Omega^5_Y|_\ell)$, but \eqref{bock} represents
it as an element of $H^1(\TTT_Y(-4)|_\ell)$, and the wedge
product becomes the contraction of 1-forms with vector fields.
Thus we are computing $\Tr(\nu,\psi)$, which is nothing but the image of ($\nu$, $\psi$) under
the Serre coupling
$$
H^1(\TTT_Y(-4)|_\ell)\times H^0(\Omega^1_Y(2)|_\ell)\lra
H^1(\OOO_\ell(-2)).
$$
On the level of \v{C}ech cocycles, the value of the coupling
on a pair $(\sum_i \phi_ie_i,\sum_j \psi_j \check e_j)$
is the coefficient of $\frac{1}{x_6x_7}$
in the expression $\sum_i\phi_i\psi_i$. This is due to the fact that the cohomology
class of the cocycle $\frac{1}{x_6x_7}$ generates $H^1(\OOO_\ell(-2))$,
and the classes of all the other cocycles of the form $\frac{1}{x_6^ix_7^j}$
($i+j=2,\ (i,j)\neq (1,1)$) are zero. We obtain that $\alp_4(m)(\xi_1,\ldots,\xi_4)$
is the coefficient of $\frac{1}{x_6x_7}$ in the product
$$
\left(\frac{x_7^2e_3+x_6x_7e_4+x_6^2e_5}{x_6^3x_7^3}\ ,
\begin{vmatrix}
a_{10}x_6+a_{11}x_7&b_{10}x_6+b_{11}x_7&c_{11}&c_{12}&c_{13}\\
     \vdots&\vdots&\vdots&\vdots\\
 a_{40}x_6+a_{41}x_7&b_{40}x_6+b_{41}x_7&c_{41}&c_{42}&c_{43}\\
 \check e_1&\check e_2&\check e_3&\check e_4&\check e_5&
\end{vmatrix}\right),
$$
which coincides with \eqref{4dets}.
\end{proof}

\begin{corollary}\label{4-form}
Let $Y$ be a smooth quartic hypersurface in $\PP^7$, and $\ell$
a line of first type in $Y$. Then $F(Y)$ is smooth at $m=[\OOO_\ell]$,
and there exists a basis ($a_0$, $a_1$, $b_0$, $b_1$, $c_1$, $c_2$, $c_3$)
of the $\OOO_{F(Y)}$-module $\Omega^1_{F(Y)}$ over an open neighborhood $U$ of $m$
such that the $4$-form $\alp_4$ defined on $U$ by Theorem \ref{closed}
is given by
$$
\alp_4=a_0\wedge b_0\wedge c_2 \wedge c_3
-a_0\wedge b_1\wedge c_1 \wedge c_3
-a_1\wedge b_0\wedge c_1 \wedge c_3
+a_1\wedge b_1\wedge c_1 \wedge c_2.
$$

\end{corollary}

\section{Nondegeneracy of exterior forms}\label{4-forms}

\textcolor{black}{
We continue to study the 4-form $\alp_4$ defined on $F(Y_4)$ for a general
quartic $Y_4$. This is the first example that goes beyond the
Beauville-Donagi case, in which $\alp_2$ is a nondegenerate 2-form, but it already
shows that the question on the ``degree of nondegeneracy''
of the exterior forms we construct is highly nontrivial.}

The classification of trilinear alternating forms in seven
complex variables goes back to Schouten (1931). There are exactly ten
orbits of $GL_7$ (including zero); in particular, as was already known to
E. Cartan, there is an open orbit of forms whose stabilizer is (up to
a finite group) a form of the exceptional group $G_2$ (see \cite{CH}).
Moreover the normal form of a generic alternating 3-form encodes
the multiplication table of the Cayley octonion algebra.

The classification of 4-forms in seven complex variables is almost equivalent
to that of 3-forms, because of the isomorphism $\wedge^4U\simeq \wedge^3U^\vee
\otimes\det(U)$ for a seven dimensional vector space $U$. More precisely,
the projective classifications are completely
equivalent, and it follows that $GL(U)$ has the same orbits in $\wedge^4U$
and in $\wedge^3U^\vee$, with isomorphic stabilizers up to finite groups.

One way to distinguish the orbits is to observe that there exists a
$GL_7$-equivariant map $S^3(\wedge^3U^\vee)\rightarrow S^2U^\vee$ defined
up to constant. Indeed, let us choose a generator $\Omega$ of $\det(U)^\vee$.
To each $\omega\in\wedge^3U^\vee$ we associate the quadratic form $q_\omega$
on $U$ defined by the formula
$$\omega\wedge (u\lrcorner\;\omega)\wedge (u\lrcorner\;\omega)
=q_\omega(u)\Omega \qquad\forall u\in U,$$
where $u\lrcorner\,\omega$ denotes the 2-form obtained by contracting
$\omega$ with $u$.
Then $q_\omega$ is non-degenerate if and only if $\omega$ belongs to the open
orbit, and this yields the classical embedding $G_2\subset SO_7$.
On the complement of this  open orbit the rank of $q_\omega$ drops to four.
More precisely, it is equal to four exactly on the codimension one orbit.

\smallskip
One can check that the orbit structure of $\wedge^4U$ is as follows,
where we denote by $\mathcal{O}_k^d(r,\rho)$ the orbit corresponding to
the 3-form denoted $f_k$ in \cite{CH}, Table 1; $d$ is the dimension
of the orbit, $r$ is the 2-rank and $\rho$ the rank.

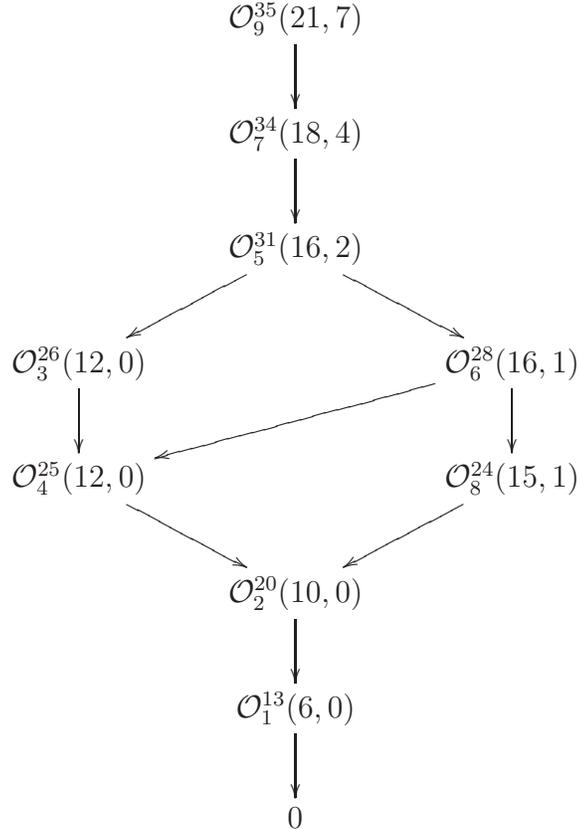
\begin{equation*}
\xymatrix
{ & \mathcal{O}_9^{35}(21,7)\ar[d] & \\
 & \mathcal{O}_7^{34}(18,4)\ar[d] & \\
 & \mathcal{O}_5^{31}(16,2)\ar[dr]\ar[dl] & \\
\mathcal{O}_3^{26}(12,0)\ar[d] & & \mathcal{O}_6^{28}(16,1)\ar[d]\ar[dll] \\
\mathcal{O}_4^{25}(12,0)\ar[dr] & & \mathcal{O}_8^{24}(15,1)\ar[dl] \\
 & \mathcal{O}_2^{20}(10,0)\ar[d] & \\
 & \mathcal{O}_1^{13}(6,0)\ar[d] & \\
 & 0 &
}\end{equation*}

Comparing the expression of our 4-form in Corollary \ref{4-form},
we see that it does not belong to the open orbit $\mathcal{O}_9^{35}$,
so it does not
define a $G_2$-structure. Nevertheless, it belongs to the codimension
one orbit $\mathcal{O}_7^{34}$, and in this  sense it can be
said ``minimally degenerate''.

\setcounter{section}{0}\renewcommand{\thesection}{\Alph{section}}
\section{Appendix. Lines on Pfaffian varieties}\label{lpv}

Let $F(\Pf_k(W^*))$ denote the Fano scheme of
lines on $\Pf_k(W^*)$. By definition it is a subscheme of
$\Gr(2,\Lambda^2W^*)$ with equations given by $\sigma_k$.
In other words, the sheaf of ideals of $F(\Pf_k(W^*))$ is the image of the map
$$
\xymatrix@1{\Lambda^{2(k-1)}W \otimes S^{-(n-k+1)}L \ar[r]^-{\sigma_k} & \CO_{\Gr(2,\Lambda^2W^*)}},
$$
where $\sigma_k$ is defined in~\eqref{sik} and $L$ is the tautological bundle on $\Gr(2,\Lambda^2W^*)$.

Consider the product $\Gr(n+k,W)\times\Gr(2,\Lambda^2W^*)$ and the subvariety
\begin{equation}\label{deftfk}
\TF_k = \{(U,L) \in \Gr(n+k,W)\times\Gr(2,\Lambda^2W^*)\ |\ L \subset \Ker(\Lambda^2W^* \to \Lambda^2U^*) \}.
\end{equation}
The projection $\pi$ of $\Gr(n+k,W)\times\Gr(2,\Lambda^2W^*)$ onto the second factor
maps $\TF_k$ into $F(\Pf_k(W^*))$. The goal of this section is to prove the following

\begin{theorem}\label{resfpf}
The map $\pi:\TF_k \to F(\Pf_k(W^*))$ is a resolution of singularities.
For $k = 1$ this resolution is crepant.
\end{theorem}

We start with the following simple observation.

\begin{lemma}\label{tfksm}
$\TF_k$ is a smooth connected algebraic subvariety of codimension $(n+k-1)(n+k)$
in $\Gr(n+k,W)\times\Gr(2,\Lambda^2W^*)$.
\end{lemma}
\begin{proof}
Consider the projection $\TF_k \to \Gr(n+k,W)$. It is clear that its fiber
over a point $U \in \Gr(n+k,W)$ is the Grassmannian $\Gr(2,\Ker(\Lambda^2W^* \to \Lambda^2U^*))$.
Therefore $\TF_k$ is smooth and connected.
This also allows to compute the codimension.
\end{proof}

The most complicated part is the surjectivity of $\pi$.

\begin{proposition}
The map $\pi:\TF_k \to F(\Pf_k(W^*))$ is surjective.
\end{proposition}
\begin{proof}
By definition $\TF_k$ is the zero locus of the canonical global section $\phi$ of the vector bundle
$\Lambda^2U^*\boxtimes L^*$ where $U$ and $L$ denote the tautological bundles
on the Grassmannians $\Gr(n+k,W)$ and $\Gr(2,\Lambda^2W^*)$ respectively.
By Lemma~\ref{tfksm} this section is regular, hence we have the Koszul resolution
\begin{equation}\label{krtf}
\dots \to \Lambda^2(\Lambda^2U \boxtimes L) \to \Lambda^2U \boxtimes L \to \CO \to \CO_{\TF_k} \to 0.
\end{equation}
Its $t$-th term equals
$$
\Lambda^t(\Lambda^2 U \boxtimes L) =
\bigoplus_{a+b=t,\ a\ge b} \Sigma^{(a,b)'}(\Lambda^2U) \boxtimes \Sigma^{(a,b)}L,
$$
where $\Sigma^\alpha$ is the Schur functor associated with the partition $\alpha$
and $\alpha'$ denotes the transposed partition of $\alpha$. Let $\pi$ be the projection
of the product $\Gr(n+k,W)\times\Gr(2,\Lambda^2W^*)$ onto the second factor.
We are going to apply the pushforward functor $\pi_*$ to the above Koszul resolution.
For this we will need the following vanishing result.

\begin{lemma}\label{cohvan}
We have $H^q(\Gr(n+k,W),\Sigma^{(a,b)'}(\Lambda^2U)) = 0$ for all $q \ge a+b - 1$ with only two exceptions:
$$
H^{n-k}(\Gr(n+k,W),\Sigma^{(n-k+1,0)'}(\Lambda^2U)) = \Lambda^{2k-2}W,
\qquad
H^{0}(\Gr(n+k,W),\Sigma^{(0,0)'}(\Lambda^2U)) = \CC.
$$
\end{lemma}

We postpone the proof of the Lemma, and now finish with the Proposition.
Note that
$$
\pi_*(\Sigma^{(a,b)'}(\Lambda^2U) \boxtimes \Sigma^{(a,b)}L) \cong
H^\bullet(\Gr(n+k,W),\Sigma^{(a,b)'}(\Lambda^2U)) \otimes \Sigma^{(a,b)}L.
$$
Therefore, the spectral sequence of hyperdirect images together with the above Lemma
shows that we have the following exact sequence
$$
\xymatrix@1{
\dots \ar[r] &
\Lambda^{2(k-1)}W \otimes S^{-(n-k+1)}L \ar[r]^-{\xi_k} &
\CO_{\Gr(2,\Lambda^2W^*)} \ar[r] &
\pi_*\CO_{\TF_k} \ar[r] & 0,}
$$
so it follows that $\pi(\TF_k)$ is the zero locus of a section $\xi_k$ of the vector bundle
$\Lambda^{2(k-1)}W^* \otimes S^{n-k+1}L^*$ on $\Gr(2,\Lambda^2W^*)$.
Moreover, since $\TF_k$ is $\GL(W)$-invariant it follows that $\xi_k$ should be $\GL(W)$-semi-invariant.
But as it was mentioned in Section~\ref{phtd}, the only $\GL(W)$-semi-invariant section there is $\sigma_k$,
so $\xi_k = \sigma_k$ and we conclude that $\pi(\TF_k) = F(\Pf_k(W^*))$. In particular, the projection $\pi$
restricts to a surjective map $\pi:\TF_k \to F(\Pf_k(W^*))$, which is precisely what is claimed in the Proposition.
\end{proof}

\begin{corollary}\label{fpfirr}
The Fano scheme $F(\Pf_k(W^*))$ is irreducible.
\end{corollary}

\begin{proof}[Proof of the Lemma~\ref{cohvan}]
Note that $\Sigma^{(a,b)'}(\Lambda^2U) \subset \Lambda^a(\Lambda^2U) \otimes \Lambda^b(\Lambda^2U)$.
The decomposition of $\Lambda^a(\Lambda^2U)$ into irreducible components is
$$\Lambda^a(\Lambda^2U)=\bigoplus_{\lambda\in W(a)}\Sigma^{d(\lambda)}U,$$
the sum being taken over the set $W(a)$ of non decreasing sequences $\lambda=(\lambda_1\ge\cdots
\lambda_c\ge c)$ such that $a=|\lambda|+\frac{c(c-1)}{2}$, and we use the notation
$$d(\lambda)=(\lambda_1,\ldots, \lambda_c,c^{\lambda_c-c+1}, (c-1)^{\lambda_{c-1}-\lambda_c},
\ldots,1^{\lambda_1-\lambda_2}).$$

Now assume that $\alpha = (\alpha_1,\alpha_2,\dots,\alpha_{n+k})$ is a Young diagram with $n+k$ rows
and even number $2m$ of boxes such that $\Sigma^\alpha U$ is a component of $\Sigma^{(a,b)'}(\Lambda^2U)$.
Then $H^q(\Gr(n+k,W),\Sigma^\alpha U) \ne 0$ only if
$$
q=p(n-k),\quad
\alpha_1\ge \dots \ge \alpha_p \ge n-k+p,\quad\text{and}\quad
p \ge \alpha_{p+1} \ge \dots \ge \alpha_{n+k}\quad
\text{for some $p$.}
$$
Moreover, $\Sigma^\alpha U$ must be a component of the tensor product of $\Sigma^{d(\lambda)}U$ and
$\Sigma^{d(\mu)}U$ for some $\lambda=(\lambda_1\ge\cdots\ge
\lambda_c\ge c) \in W(a)$ and some $\mu=(\mu_1\ge\cdots\ge
\mu_d\ge d)\in W(b)$.

If $p = 0$ then $\alpha = 0$, so we get the case $H^{0}(\Gr(n+k,W),\Sigma^{(0,0)'}(\Lambda^2U)) = \CC$.

Assume that $p\ge 1$. Recall that by the Littlewood-Richardson rule, the diagram of $\alpha$ is obtained
by adding $2b$ boxes to the diagram of $d(\lambda)$ (with certain restrictions imposed by $d(\mu)$).
In particular $\alpha_c\ge d(\lambda)_c\ge c$, and therefore $p\ge c$ (and symmetrically $p\ge d$).
Because of the special form of the partition $d(\lambda)$, this implies that
$$\begin{array}{rcl}
d(\lambda)_{p+1}+\cdots +d(\lambda)_{n+k} & = & \frac{1}{2}|d(\lambda)|-\frac{c(c+1)}{2}-(
d(\lambda)_{c+2}+\cdots +d(\lambda)_{p}) \\
 & \ge & a-\frac{c(c+1)}{2}-c(p-c-1)=a+\frac{c(c+1)}{2}-cp.
\end{array}$$
And of course we have a similar estimate for $d(\lambda)$.

Now, we observe that following the Littlewood-Richardson rule, the diagram of $\alpha$ is obtained from that
of $d(\lambda)$ by adding some numbered boxes, with for each $i$, exactly $d(\mu)_i$ boxes numbered $i$.
Moreover the $j$-th line can only contain
boxes numbered by some $i\le j$, so all the boxes numbered $k$ must appear on line $k$ or below.
This implies that
$$\begin{array}{rcl}
\alpha_{p+1}+\cdots +\alpha_{n+k}  & \ge &
d(\lambda)_{p+1}+\cdots +d(\lambda)_{n+k}+d(\mu)_{p+1}+\cdots +d(\mu)_{n+k} \\
 & \ge & a+\frac{c(c+1)}{2}-cp+b+\frac{d(d+1)}{2}-dp.
\end{array}$$
But since $\alpha_1\ge \dots \ge \alpha_p \ge n-k+p$, we deduce that
$$\begin{array}{rcl}
a+b & \ge &  p(n-k+p)+\frac{c(c+1)}{2}-cp+\frac{d(d+1)}{2}-dp\\
 &= & q+\frac{(p-c)^2}{2}+\frac{(p-d)^2}{2}+\frac{c+d}{2}\ge q+p.
\end{array}$$
Under the hypothesis that $q\ge a+b-1$, we deduce that $p=c=d=1$,
and the only possibility is
$\alpha = (n-k+1,1,\dots,1,0,\dots,0)$ (the number of ones being $n-k+1$).
Then $H^{n-k}(\Gr(n+k,W),\Sigma^\alpha U) = \Lambda^{2k-2}W$, which gives the second case, since
$\Sigma^\alpha U$ is a component of $\Sigma^{(n-k+1,0)'}(\Lambda^2U)$ (with multiplicity one),
but of no other $\Sigma^{(a,b)'}(\Lambda^2U)$.
\end{proof}

\begin{proposition}\label{pibir}
The map $\pi:\TF_k \to F(\Pf_k(W^*))$ is birational.
\end{proposition}
\begin{proof}
Since $\pi$ is surjective and proper, and $F(\Pf_k(W^*))$ is irreducible,
it suffices to check that it is an isomorphism over an open subset of $F(\Pf_k(W^*))$.
Let $F_k^0 \subset F(\Pf_k(W^*))$ be the open subset consisting of lines which
do not intersect the locus of skew-forms of rank $\le 2n-2k-2$. Let $L$ be such a line.
Note that if $U \subset W$ is a subspace of dimension $n+k$ such that $\phi_{L,U} = 0$,
then for each point in $L$ the space $U$ contains a $2k$-dimensional subspace lying
in the kernel of the corresponding skew-form. But if all skew-forms in $L$ have rank
$2n-2k$ then $U$ contains all their kernels. Thus it suffices to check that for such $L$
the linear hull of the kernels of all forms in $L$ has dimension $n+k$. For this we consider
the following exact sequence on $\PP(L)$:
$$
0 \to K(-1) \to W\otimes\CO_{\PP(L)}(-1) \to W^*\otimes\CO_{\PP(L)} \to K^* \to 0,
$$
where the middle map is the evaluation of a skew-form on a vector and $K$ is the bundle of kernels
of skew-forms. Then the codimension of the linear
hull of kernels coincides with $H^0$ of the image $I$ of the middle map.
Computing the determinant of this exact sequence we see that $\deg K^* = n - k$.
Since $K^*$ is a vector bundle of rank $2k$ on $\PP(L) = \PP^1$ generated by global sections
we conclude by Riemann--Roch that $\dim H^0(\PP(L),K^*) = n + k$. On the other hand, from the above
exact sequence it follows that the map $H^0(\PP(L),W^*\otimes\CO_{\PP(L)}) \to H^0(\PP(L),K^*)$
is surjective, hence $\dim H^0(\PP(L),I) = 2n - (n+k) = n-k$. So the codimension of the linear hull
of kernels is $n-k$, hence the dimension is $n+k$ as required.
\end{proof}

The next goal is to investigate where the map $\pi$ is not an isomorphism.
For this we need the following.

\begin{lemma}
Let $L \in F(\Pf_k(W^*))$ be a line on $\Pf_k(W^*)$.
If $L \cap \Pf_{k+1}(W^*) \ne \emptyset$ then $\pi^{-1}(L)$ contains a line.
\end{lemma}
\begin{proof}
Let $U \subset W$ be a subspace of dimension $n+k$ isotropic for $L$
and assume that the form $\lambda$ corresponding to a point of $L$
has rank strictly less than $2n-2k$. Then there exists a vector subspace
$U' \subset W$ of dimension $n+k+1$ containing $U$ and isotropic for $\lambda$. Let $U' = U \oplus \CC w$.
Let $\lambda'$ be another skew-form in $L$. Then $\lambda'(-,w)$ is a linear form on $U$.
Let $U''$ be its kernel. Then it is clear that each vector subspace of dimension
$n+k$ in $U'$ containing $U''$ is isotropic both for $\lambda$ and $\lambda'$.
These subspaces form a line in $\Gr(n+k,W)$ which is contained in $\pi^{-1}(L)$.
\end{proof}

By Lemma~\ref{tfksm} variety $\TF_k$ is smooth and by Proposition~\ref{pibir}
the map $\pi:\TF_k \to F(\Pf_k(W^*))$ is birational, hence $\TF_k$ is a resolution of $F(\Pf_k(W^*))$.
The final statement of the Theorem is given by the following

\begin{lemma}
The resolution $\pi:\TF_1 \to F(\Pf(W^*))$ is crepant.
\end{lemma}
\begin{proof}
Compute the canonical class of $\TF_1$. Recall that $\TF_1$ is the zero locus
of a regular section of $\Lambda^2U^*\boxtimes L^*$ on $\Gr(n+1,W)\times\Gr(2,\Lambda^2W^*)$.
We have
$$
\begin{array}{l}
\omega_{\Gr(n+1,W)\times\Gr(2,\Lambda^2W^*)} = \CO(-2n,-n(2n-1)),\\
\det(\Lambda^2U^*\boxtimes L^*) =
\det(\Lambda^2U^*)^{\otimes 2} \otimes (\det L^*)^{\otimes n(n+1)/2} =
\CO(2n,n(n+1)/2),
\end{array}
$$
so by adjunction we get $\omega_{\TF_1} \cong \CO(0,-3n(n-1)/2)$.
Thus the canonical class of $\TF_1$ is a pullback from $F(\Pf(W^*)) \subset \Gr(2,\Lambda^2W^*)$,
so $\pi:\TF_1 \to F(\Pf(W^*))$ is crepant.
\end{proof}

\section{Appendix: proof of Proposition \ref{Hodge_num}}\label{proof_betti_numbers}

\textcolor{black}{
The Hodge numbers of $Y$ are easily obtained via the Griffiths residue
isomorphism. So we are turning to the calculation of the Hodge numbers of $F$
obtained by an application of the Borel-Weil theorem to some concrete homogeneous bundles.
Though this calculation is somewhat lengthy and takes two pages, we think it is
worthwhile to include it, because it should be very useful
to those algebraic geometers that have no experience of practical application of
representation theory. Moreover, it does not only give Hodge numbers,
but represents the respective Hodge groups as irreducible representations of
certain linear groups, which may help to find
isomorphic pieces of Hodge structures in different varieties. This is
important in view of the problem of searching
for {\em special} $Y_n$'s in the sense of the definition given in the Introduction (or
in \cite{Ha} when $n=3$).
}

To compute the Hodge numbers of $F$, we will use the fact that $F$ is the zero-locus
in the Grassmannian $G=\Gr(2,8)$ of a general section of the vector bundle
$S^4T^*$, where $T$ is the tautological rank two bundle on
$G$. This bundle has rank five and is generated by global sections.
We can therefore use the conormal sequence
$$0\rightarrow S^4T_{|F}\rightarrow\Omega^1_{G|F} \rightarrow
\Omega^1_{F}\rightarrow 0,$$
as well as the Koszul complex
$$0\rightarrow \wedge^5(S^4T)\rightarrow \wedge^4(S^4T)\rightarrow
\wedge^3(S^4T)\rightarrow \wedge^2(S^4T)
\rightarrow S^4T\rightarrow\mathcal{O}_G\rightarrow\mathcal{O}_F\rightarrow 0$$
in order to compute the cohomology of the restriction to $F$ of
vector bundles  on $G$. The exterior powers of $S^4T$ are readily computed
in terms of Schur powers (note that since $T$ has rank two, we simply have
$\Sigma^{a,b}T=S^{a-b}T\otimes\mathcal{O}(-b)$):
$$\begin{array}{lcl}
\wedge^2(S^4T) & = & \Sigma^{7,1}T\oplus \Sigma^{5,3}T \\
\wedge^3(S^4T) & = & \Sigma^{9,3}T\oplus \Sigma^{7,5}T \\
\wedge^4(S^4T) & = & \Sigma^{10,6}T \\
\wedge^5(S^4T) & = & \Sigma^{10,10}T
\end{array}$$
In order to compute the cohomology of $S^4T$ and $\Omega^1_{G}$
restricted to $F$, we will twist the Koszul resolution of $\mathcal{O}_F$ by
these bundles and use Bott's theorem on $G$. Recall that if $Q$ denotes
the rank $6$ quotient vector bundle on $G$, we have $\Omega^1_{G}\simeq
Q^*\otimes T$. Bott's theorem computes the cohomology of any tensor
product of Schur powers of $Q$ and $T$ as follows. Let $\alpha$ and $\beta$
be two non-increasing sequences of relative integers, of respective lengths
$6$ and $2$. Let $\rho=(8,7,6,5,4,3,2,1)$ and consider the sequence
$(\alpha,\beta)+\rho$. Call it {\it regular} if its entries are pairwise
distinct. In that case there is a unique permutation $w$ such that
$w((\alpha,\beta)+\rho)$ is strictly decreasing. Then $\lambda =
w((\alpha,\beta)+\rho)-\rho$ is non-increasing. Bott's theorem asserts,
if $V$ is the ambient eight-dimensional space,  that
$$H^q(G,\Sigma^{\alpha}Q\otimes \Sigma^{\beta}T)=\Sigma^{\lambda}V$$
if $(\alpha,\beta)+\rho$ is regular and $q=\ell (w)$, and zero otherwise.
For example, if $\beta=(a,b)$, we get that $H^q(G,\Sigma^{\beta}T)\neq 0$
if and only if either
\begin{eqnarray*}
 q=0, &\quad &0\ge a\ge b; \\
 q=6, &\quad &a\ge 7, 1\ge b; \\
 q=12, &\quad &a\ge b\ge 8.
\end{eqnarray*}
We can easily deduce the cohomology groups of $\mathcal{O}_F$.
Indeed, $H^q(G,\wedge^i(S^4T))$ is non zero only for $(i,q)=(0,0),
(2,6),(5,12)$, which implies that $H^q(F,\mathcal{O}_F)$ is non
zero exactly for $q=0,4,7$. In particular
$$H^4(F,\mathcal{O}_F)\simeq H^6(G,\wedge^2(S^4T))\simeq\CC.$$

Let us turn to the cohomology groups of the cotangent bundle of $F$.
Bott's theorem
implies that the only non zero groups among the $H^q(G,S^4T\otimes
\wedge^i(S^4T))$ appear in bidegree $(i,q)=(1,6),(2,6),(4,12),(5,12)$.
In particular $q-i$ is always bigger than three and we can deduce that
$H^q(F,S^4T_{|F})=0$ for $q\le 3$. Moreover, Bott's theorem gives
$$\begin{array}{lcl}
H^6(G,S^4T\otimes S^4T) &= & End(V),\\
H^6(G,S^4T\otimes \wedge^2(S^4T)) &= &S^4V,
\end{array}$$
and we deduce an exact sequence
$$0\rightarrow H^4(F,S^4T_{|F})\rightarrow S^4V\rightarrow End(V)\rightarrow H^5(F,S^4T_{|F})\rightarrow 0.$$
The middle map $S^4V\stackrel{\psi_f}{\rightarrow} End(V)$
must be dual to the map $End(V)\stackrel
{\phi_f}{\rightarrow} S^4V^*$ mapping $u$ to $u(f)$, where $f$ denotes an equation
of $Y$. Indeed, we can do the same computation in family, with a variable
$f$, and use the $GL(V)$-equivariance to ensure that the map $\phi_f$
depends linearly on $f$. And there is, up to scalar, a unique
equivariant map from $S^4V^*$ to $Hom(End(V),S^4V^*)$.
We can conclude that for $f$ general, $\psi_f$ is surjective, and that its
kernel  $H^4(F,S^4T_{|F})\simeq H^{4,2}(Y)^*$. Indeed, this is just a
reformulation of the Griffiths isomorphism.

There remains to compute the cohomology groups of $\Omega^1_{G|F}$.
Applying Bott's theorem as above, one checks that
$H^q(G,\Omega^1_G\otimes \wedge^i(S^4T))$ is non zero only
for $(i,q)=(0,1),(2,7),(5,12)$, which implies that $H^q(F,\Omega^1_{G|F})$
is non zero only for $q=1,5,7$. But then we can deduce from the conormal
exact sequence that
$$H^{1,3}(F)\simeq H^4(F,S^4T_{|F})\simeq H^{2,4}(Y).$$

Finally we need to compute $H^{2,2}(F)$. Let us denote by $K$ the kernel
of the natural map $\Omega^2_{G|F}\rightarrow \Omega^2_{F}$. Bott's theorem
implies that $H^q(G,\Omega^2_G\otimes \wedge^i(S^4T))$ is non zero
only for $(i,q)=(0,2),(2,7),(2,8),(3,8),(4,12),(5,12)$. In particular
the only non zero group $H^q(F,\Omega^2_{G|F})$ for $q\le 4$ is
$$H^2(F,\Omega^2_{G|F})\simeq H^2(G,\Omega^2_G)=\CC^2,$$
and we get an exact sequence
$$H^2(F,K)\rightarrow H^{2,2}(G)\rightarrow H^{2,2}(F)\rightarrow H^3(F,K)\rightarrow 0.$$
In order to compute the cohomology groups of $K$, consider the exact sequence
$$0\rightarrow S^2(S^4T)_{|F}\rightarrow S^4T\otimes \Omega^1_{G|F}\rightarrow K\rightarrow 0.$$
Applying Bott's theorem once again, we check that $H^q(G,S^4T\otimes
\Omega^1_G\otimes \wedge^i(S^4T))$ is non zero only for $(i,q)=(1,6),
(1,7),(2,7),(3,12),(4,12)$, which implies that $H^q(F,S^4T_{|F}\otimes
\Omega^1_{G|F})=0$ for $q\le 4$. Therefore $H^q(F,K)\simeq H^{q+1}(F,
S^2(S^4T)_{|F})$ for $q\le 3$. To compute the latter, we consider
the groups $H^q(G,S^2(S^4T)\otimes \wedge^i(S^4T))$ and check by Bott's
theorem that they are non zero exactly for $(i,q)=(0,6),(1,6),(2,6)$
or $q=12$. We have
$$\begin{array}{rcl}
H^6(G,S^2(S^4T)) & = & \Sigma^{21111110}V\simeq sl(V), \\
H^6(G,S^2(S^4T)\otimes S^4T) & = & \Sigma^{61111110}V\oplus \Sigma^{51111111}V
\simeq S^5V\otimes V^*,\\
H^6(G,S^2(S^4T)\otimes \wedge^2(S^4T)) & = & \Sigma^{91111111}V\simeq S^8V.
\end{array}$$
Hence the exact sequence
$$0\rightarrow H^3(F,K)\rightarrow S^8V\stackrel{\alpha_f}{\rightarrow} S^5V\otimes V^*.$$
As above we can argue that the map $\alpha_f$ must be dual to the natural
map $\beta_f$ obtained as the composition
$$S^5V\otimes V^*\rightarrow S^5V\otimes S^3V\rightarrow S^8V$$
deduced from the map $V^*\rightarrow S^3V$ given by the differential of $f$.
Otherwise said, we are just multiplying quintic polynomials with the
derivatives of $f$, so that Griffiths' residue theorem tells us precisely
that the cokernel of $\beta_f$ is isomorphic to $H^{3,3}(Y)_{prim}$.
Since this cokernel is dual to the kernel of $\alpha_f$, we finally get
$$H^{2,2}(F)\simeq H^{2,2}(G)\oplus H^3(F,K)\simeq
H^{2,2}(G)\oplus H^{3,3}(Y)_{prim}.$$

\end{document}